\def\technicalreport{1}

\documentclass[10pt]{article}

\RequirePackage{amsmath}
\RequirePackage{amsfonts}
\RequirePackage{amssymb}
\RequirePackage{hyperref}
\RequirePackage{multirow}
\RequirePackage{algorithm}
\RequirePackage{algpseudocode}
\RequirePackage{graphicx}
\RequirePackage{xcolor}
\RequirePackage{amsthm}
\RequirePackage{amsmath}
\RequirePackage{amsfonts}
\RequirePackage{amssymb}
\RequirePackage{hyperref}
\RequirePackage{multirow}
\RequirePackage{algorithm}
\RequirePackage{algpseudocode}
\RequirePackage{graphicx}
\RequirePackage{xcolor}
\allowdisplaybreaks

\usepackage{dpr_fec,ifthen,subcaption,dsfont}

\usepackage{isetechreport}
\def\reportyear{19T}
\def\reportno{025}
\def\revisionno{0}
\def\originaldate{\today}
\coraltrue
\cvcrfalse

\begin{document}

\title{A Fully Stochastic Second-Order Trust Region Method}

\author{Frank E.~Curtis\thanks{E-mail: frank.e.curtis@lehigh.edu}}
\author{Rui Shi\thanks{E-mail: rus415@lehigh.edu}}
\affil{Department of Industrial and Systems Engineering, Lehigh University}
\titlepage

\maketitle

\begin{abstract}
  A stochastic second-order trust region method is proposed, which can be viewed as a second-order extension of the \emph{trust-region-ish} (TRish) algorithm proposed by Curtis et al.~[INFORMS J.~Optim.~1(3) 200--220, 2019].  In each iteration, a search direction is computed by (approximately) solving a trust region subproblem defined by stochastic gradient and Hessian estimates.  The algorithm has convergence guarantees for stochastic minimization in the fully stochastic regime, meaning that guarantees hold when each stochastic gradient is required merely to be an unbiased estimate of the true gradient with bounded variance and when the stochastic Hessian estimates are bounded uniformly in norm.  The algorithm is also equipped with a worst-case complexity guarantee in the nearly deterministic regime, i.e., when the stochastic gradient and Hessian estimates are very close in expectation to the true gradients and Hessians.  The results of numerical experiments for training convolutional neural networks for image classification and training a recurrent neural network for time series forecasting are presented.  These results show that the algorithm can outperform a stochastic gradient approach and the first-order TRish algorithm in practice.
\end{abstract}

\newcommand{\finf}{f_{\rm inf}}

\section{Introduction}

For many years, the foundational approach for solving stochastic optimization problems has been the stochastic gradient method (\cite{RobbMonr51}), hereafter referred to as~SG.  However, despite its many theoretical and practical advantages, there remain some shortcomings in the use of SG for solving many stochastic optimization problems, including many arising in machine learning and signal processing, areas in which SG and its variants are very popular.  For example, one disadvantage of~SG and many variants of it (see \S\ref{sec.literature_review}) is that the variance of the step taken by the algorithm in each iteration is proportional to the variance of the stochastic gradient estimate, which can be large.  In the \emph{fully stochastic} regime, i.e., when the variances of the stochastic gradient estimates are merely bounded by some (large) constant, SG can take a large step even though the norm of the true gradient may be relatively small in norm.

In~\cite{CurtScheShi19}, a first-order stochastic optimization algorithm is proposed that is designed to mitigate the effects of large variances of the stochastic gradient estimates.  Based on a trust region methodology, this \emph{trust-region-ish} algorithm, known as TRish, uses a careful step normalization procedure in order to attain theoretical convergence properties on par with those of SG, but in such a way that the empirical performance can be better than that of SG.  The results of experiments on logistic regression and deep neural network training problems demonstrates that the empirical performance of TRish can be significantly better than that of a traditional SG approach.  In particular, TRish is able to reach better solutions more quickly, and in a more stable manner, meaning that the quality of the solution estimate does not vary wildly from one iteration to the next.

In this paper, we extend the TRish methodology to allow for the use of stochastic second-order information, in the form of stochastic Hessian estimates that are incorporated in the trust region subproblems.  The resulting algorithm, which we continue to refer to as TRish, is shown to have good convergence properties in a wide range of settings.  In particular, in the fully stochastic regime and with a very loose requirement on the accuracy with which the trust region subproblems are solved, we show that the algorithm achieves convergence properties on par with those of TRish.  Admittedly, this is done with assumptions that impose stricter requirements on the stepsizes employed in the algorithm, but the results are still non-trivial to obtain, and the theoretical analysis in this paper requires different techniques than those employed in \cite{CurtScheShi19}.  We also include some theoretical guarantees that are stronger than have been presented for the first-order variant of TRish.  On the other end of the theoretical spectrum, we show that when the stochastic gradient and Hessian estimates are very close in expectation to the true values, and when the subproblems are solved exactly, TRish offers a worst-case complexity property that is similar to that offered by a determinstic second-order trust region method.

As has been the motivation for other authors considering second-order extensions of stochastic optimization algorithms, a main motivation for our work is to design an algorithm that can ideally inherit the benefits of Newton-trust-region methods for minimization, such as their scale invariance, ability to employ problem-independent stepsizes near a solution (as opposed to first-order methods, which require problem-dependent stepsizes, typically related to the Lipschitz constant for the gradient of the objective), ability to handle nonconvexity and avoid saddle points without extra computational procedures, and asymptotic fast rate of convergence.  These properties cannot fully be attained in the stochastic regime, but our numerical experiments demonstrate that the TRish methodology can benefit from the use of stochastic second-order derivative information in practice.  The results that we present in this paper are for training convolutional neural networks (CNNs) for image classification, and for training a recurrent neural network (RNN) for time series forecasting.  Our results suggest that TRish can be an effective approach for stochastic and finite-sum minimization over broad classes of challenging problems.

\section{Literature Review}\label{sec.literature_review}

The literature on SG, a stochastic first-order method, is extensive.  For a few examples of papers with analyses of SG and variants of it, see~\cite{BottCurtNoce18,rakhlin2012making,AgarBott15,ByrdChinNoceWu12,Chun54,DanLan15,FrieSchm12,NemiJudiLanShap09,GhadLan13,RobbMonr51}, and \cite{RobbSieg71}.

Stochastic second-order methods, which can be classified as methods that compute each step by (approximately) minimizing a quadratic model of the objective function, have received less attention in the literature.  That said, many types of methods have been proposed, analyzed, and tested.  Overall, one may characterize stochastic second-order methods into four categories (see \cite{BottCurtNoce18}): stochastic Newton methods, stochastic quasi-Newton methods, natural gradient methods, and diagonal-scaling methods.

Stochastic Newton methods, like the deterministic Newton method for minimization, compute each step by approximately minimizing a quadratic model of the form $g_k^Ts + \thalf s^TH_ks$ over $s \in \R{n}$, where $g_k$ is a stochastic gradient estimate and $H_k$ is a stochastic Hessian estimate.  For practical purposes, such an approach would typically use an iterative method such as the conjugate gradient (CG) algorithm to minimize this quadratic function approximately.  In this manner, one need not form nor factor the matrix $H_k$; instead, one need only perform matrix-vector products with $H_k$, which can be done with back propagation.  (In nonconvex settings, a regularization term might also be added if $H_k$ might not be positive definite, or one might terminate CG once negative curvature is detected, as in the standard Steihaug-CG routine (\cite{Stei83}).)  For examples of stochastic Newton methods in the literature, see \cite{anbar1978stochastic,byrd2016stochastic,BottCurtNoce18,bollapragada2018exact}, and \cite{CurtRobi19}.

Stochastic quasi-Newton methods borrow the idea from the deterministic optimization literature that, instead of employing second-derivative information, one could derive (inverse) Hessian approximations by observing differences in gradients from one iteration to the next.  In the stochastic regime, such an approach needs to have safeguards to account for the fact that the gradients are only estimated in each iteration.  For examples of stochastic quasi-Newton methods, see \cite{schraudolph2007stochastic,byrd2016stochastic,BottCurtNoce18,Curt16}, and \cite{WangMaGoldLiu17}.

Motivated by insights from information geometry, the idea of the natural gradient method is to employ the Fischer information matrix in place of the Hessian when computing a search direction.  Due to various simplications that are required to derive a practical algorithm, such an approach reduces to a type of (generalized) Gauss-Newton algorithm.  For further information in the natural gradient method and related ideas, see  \cite{martens2014new,desjardins2015natural,grosse2016kronecker}, and \cite{zhang2019fast}.

Diagonal-scaling methods, wherein each step can be expressed as a diagonal scaling matrix times the negative stochastic gradient, are not always classified as second-order methods.  However, we argue, as some others do, that these methods should be viewed in this light, and one can argue that the good performance of such methods in practice is because the algorithms are emulating second-order-type properties.  A few popular diagonal scaling methods are RMSprop (\cite{TielHint12}), Adagrad (\cite{DuchHazaSing11}), and Adam (\cite{kingma2014adam} and \cite{ReddKaleKuma18}).

Finally, it should be mentioned that some other types of stochastic trust region methods have appeared in the literature.  For a couple of examples, see \cite{chen2015stochastic} and \cite{cartis2018global}.  However, these approaches are distinct from the TRish methodology, e.g., due to their need to impose stronger requirements on the stochastic gradient and Hessian estimates to achieve their desired convergence rate guarantees.  We also direct the reader to the use of trust regions in reinforcement learning; see, e.g., \cite{schulman2015trust} and \cite{kurutach2018model}.  This setting is distinct from the one considered in this paper, but these works provide further evidence of how optimization algorithms based on trust region ideas can be effective in various settings.

\section{Algorithm}

In this section, we formally present our problem of interest, introduce relevant notation and terminology, and present our proposed algorithm.

\subsection{Problem Description}

The algorithm that we propose is designed to solve stochastic optimization problems; in particular, it is designed to minimize an objective function $f: \R{n} \rightarrow \R{}$ that is defined by an expectation of a function $F: \R{n} \times \Xi \rightarrow \R{}$ that depends on a random variable $\xi$, as in
\bequation\label{prob.f}
  \min_{x\in\R{n}}\ f(x),\ \ \text{where}\ \ f(x) = \E_\xi[F(x,\xi)].
\eequation
Here, $\E_\xi[\cdot]$ denotes expectation with respect to the distribution of~$\xi$.  A related type of problem is one obtained by taking a stochastic average approximation (SAA) of \eqref{prob.f}.  This leads to a finite-sum objective of the form $f(x) = \frac1N \sum_{i=1}^N f_i(x)$, where $f_i := F(\cdot,\xi_i)$ with~$\xi_i$ denoting a realization of the random variable~$\xi$.  Our algorithm automatically extends to this setting---whether or not the function arises from an SAA of \eqref{prob.f}---where in place of the distribution of $\xi$ one can consider a discrete uniform distribution over $\{1,\dots,N\}$.

The algorithm that we propose makes use of stochastic gradient and stochastic Hessian estimates that, at an algorithm iterate $x_k \in \R{n}$, are intended to approximate $\nabla f(x_k)$ and $\nabla^2 f(x_k)$, respectively.  These can be understood as follows.  First, in the context of \eqref{prob.f}, a stochastic gradient estimate may be computed as $g_k = \nabla_x F(x_k,\xi_k)$, where $\xi_k$ is a realization of $\xi$.  On the other hand, in the context of minimizing a finite sum, one may consider $g_k = \nabla_x f_{i_k}(x_k)$, where $i_k$ has been generated from a discrete uniform distribution over the index set $\{1,\dots,N\}$.  In either setting, $g_k$ could instead represent an average of such quantities and still be thought of as a stochastic gradient estimate.  In this case, $g_k$ is commonly referred to as a \emph{mini-batch} estimate.  Specifically, for \eqref{prob.f} one may consider $g_k = \frac{1}{|\Scal_k|} \sum_{j\in\Scal_k} \nabla_x F(x_k,\xi_{k,j})$ and for the finite-sum setting one may consider the estimate $g_k = \frac{1}{|\Scal_k|} \sum_{j\in\Scal_k} \nabla_x f_{i_{k,j}}(x_k)$.  In the statement of our algorithm, we capture all of these possibilities by writing $g_k \approx \nabla f(x_k)$.

For the stochastic Hessian estimates employed in our algorithm, we write $H_k \approx \nabla^2 f(x_k)$, but in this context the meaning of ``estimate'' is meant much more loosely.  Indeed, in the context of computing $g_k$, the possibilities in the previous paragraph make sense since our analysis requires that $g_k$ be an unbiased estimator of $\nabla f(x_k)$.  However, our assumption on $H_k$ can be much less restrictive.  While one might choose in the context of \eqref{prob.f} to define $H_k = \nabla_{xx}^2 F(x_k,\xi_k^H)$ for some realization $\xi_k^H$ of $\xi$ (or with a mini-batch), most of our analysis merely requires that $\{H_k\}$ is uniformly bounded in norm.

\subsection{Algorithm Description}

Our algorithm is stated below as \ref{alg.trish}.  Similar to the first-order version in \cite{CurtScheShi19}, each iteration involves solving a trust region subproblem involving stochastic derivative estimates.  Importantly, for much of our analysis, the algorithm merely requires that each subproblem is solved such that \emph{Cauchy decrease} is achieved.  This only requires that the solution vector~$s_k$ is feasible and yields a value for the subproblem objective that is at least as good as that offered by the \emph{Cauchy point}, which is the minimizer of the subproblem objective along its steepest descent direction from the origin (subject to the trust region constraint); see \cite{ConnGoulToin00} or \cite{NoceWrig06}.  If $H_k = 0$ for all $k \in \N{}$, then the algorithm reduces to that in \cite{CurtScheShi19}. However, clearly, the algorithm presented here offers much more computational flexibility.

\begin{algorithm}[ht]
  \renewcommand{\thealgorithm}{\texttt{TRish}}
  \caption{: (Second-Order) Trust-Region-ish Algorithm}
  \label{alg.trish}
  \begin{algorithmic}[1]\label{alg.str}
    \State Choose an initial iterate $x_1 \in \R{n}$ and positive stepsizes $\{\alpha_k\}$.
    \State Choose positive parameters $\{\gamma_{1,k}\}$ and $\{\gamma_{2,k}\}$ such that $0 < \gamma_{2,k} \leq \gamma_{1,k} < \infty$ for all $k \in \N{}$.
    \For{\textbf{all} $k \in \N{}$}
      \State Generate a stochastic gradient $g_k \approx \nabla f(x_k)$. 
      \State Compute $s_k$ yielding at least Cauchy decrease for the subproblem
      \bequation\label{prob.sub}
        \min_{s \in \R{n}}\ g_k^Ts + \half s^TH_ks\ \ \st\ \ \|s\|_2 \leq \Delta_k
      \eequation
      \State using matrix-vector products with a stochastic Hessian $H_k \approx \nabla^2 f(x_k)$, where
      \bequation\label{eq.Delta}
        \Delta_k \gets \bcases \gamma_{1,k} \alpha_k \|g_k\|_2 & \text{if $\|g_k\|_2 \in \left[0,\frac{1}{\gamma_{1,k}}\)$} \\ \alpha_k & \text{if $\|g_k\|_2 \in \left[\frac{1}{\gamma_{1,k}},\frac{1}{\gamma_{2,k}}\right]$} \\ \gamma_{2,k} \alpha_k \|g_k\|_2 & \text{if $\|g_k\|_2 \in \(\frac{1}{\gamma_{2,k}},\infty\right]$.} \ecases
      \eequation
      \State Set $x_{k+1} \gets x_k + s_k$.
    \EndFor
  \end{algorithmic}
\end{algorithm}

Further motivation for the scheme for choosing the trust region radii, namely, \eqref{eq.Delta}, can be found in \cite{CurtScheShi19}.  In short, if one were merely to choose $\Delta_k = \alpha_k$ for all $k \in \N{}$ so that the steplength is normalized in all iterations, then one might not have a convergent algorithm; it is possible that the algorithm would compute a direction that is one of expected ascent.  An example showing this possibility is shown as \cite[Ex.~1]{CurtScheShi19}.  Hence, \eqref{eq.Delta} embodies a \emph{careful} step normalization strategy that might choose $\Delta_k = \alpha_k$, but otherwise uses a nonlinear stepsize control scheme to adjust the steplength.  The specific formulas for the radii in \eqref{eq.Delta} ensure that (in the case $H_k = 0$) the steplength $\|x_{k+1} - x_k\|_2$ is a continuous function of $\|g_k\|_2$; see \cite[Fig.~1]{CurtScheShi19}.

\section{Convergence Analysis}

We prove convergence results for \ref{alg.trish} under various settings.  We begin by proving fundamental lemmas under basic sets of assumptions.  These results illuminate the critical features of the algorithm that lead to all convergence guarantees.  We present these guarantees first for the case of nonconvex $f$ and different stepsize and parameter choices, then for the case of $f$ satisfying the well-known Polyak-\L{}ojasiewicz (PL) condition, of which strongly convex functions are a special case.  Again, these results are presented for a few stepsize and parameter choices.  As \ref{alg.trish} generalizes the first-order algorithm proposed in \cite{CurtScheShi19}, the convergence theorems proved in this section essentially generalize those results proved for the first-order algorithm.  However, the proofs presented here require different approaches due to the influence of $\{H_k\}$ on the subproblems.

For convenience throughout our analysis, we denote for all $k \in \N{}$ the following cases, which clearly correspond to the different cases for the trust region radius $\Delta_k$ in \eqref{eq.Delta}:
\begin{align}
  \|g_k\| &\in \left[0,\tfrac{1}{\gamma_{1,k}}\), \tag{Case 1}\label{case.1} \\
  \|g_k\| &\in \left[\tfrac{1}{\gamma_{1,k}},\tfrac{1}{\gamma_{2,k}}\right], \tag{Case 2}\label{case.2} \\
  \text{or}\ \ \|g_k\| &\in \(\tfrac{1}{\gamma_{2,k}},\infty\). \tag{Case 3}\label{case.3}
\end{align}
Also, for shorthand, we use $\E_k[\cdot]$ to denote expectation of a random variable conditioned on the event that the algorithm has reached the iterate $x_k$ and generates the stochastic gradient estimate $g_k$ and the stochastic Hessian estimate $H_k$; i.e.,
\bequationNN
  \E_k[\cdot] \equiv \E[\ \cdot\ |\ \text{the current iterate and derivative estimates are $x_k$, $g_k$, and $H_k$}].
\eequationNN

We make the following assumptions throughout our analysis.  These assumptions are essentially the same as the basic assumptions from \cite{CurtScheShi19}, except that we add the assumption that $f$ is twice continuously differentiable, which is a reasonable assumption to add in the context of a second-order-type algorithm.

\bassumption\label{ass.f}
  The objective function $f : \R{n} \to \R{}$ is twice continuously differentiable and bounded below by a scalar $\finf := \inf_{x\in\R{n}} f(x) \in \R{}$.  In addition, the gradient function $\nabla f : \R{n} \to \R{n}$ is Lipschitz continuous with constant $L_g \in \R{}_{>0}$ $($i.e., $f$ is \emph{$L_g$-smooth}$)$.
\eassumption

It is known (see, e.g., \cite[Lemma~1.2.2]{Nest04}), that under Assumption~\ref{ass.f} one has
\bequation\label{eq.Hessian_bounded}
  \|\nabla^2 f(x)\|_2 \leq L_g\ \ \text{for all}\ \ x \in \R{n}.
\eequation

\bassumption\label{ass.gH}
  For all $k \in \N{}$, the stochastic gradient estimate $g_k$ is an unbiased estimator of the gradient $\nabla f(x_k)$ in the sense that $\E_k[g_k] = \nabla f(x_k)$.
\eassumption 

Under Assumption~\ref{ass.gH}, one finds for all $k \in \N{}$ that
\bequation\label{eq.unbiased_consequence}
  \baligned
    \E_k[\|\nabla f(x_k) - g_k\|_2^2]
      &= \E_k[\|\nabla f(x_k)\|_2^2 - 2\nabla f(x_k)^Tg_k + \|g_k\|_2^2] \\
      &= -\|\nabla f(x_k)\|_2^2 + \E_k[\|g_k\|_2^2].
  \ealigned
\eequation
  
\subsection{Fundamental Lemmas}\label{sec.fundamental}

Our first lemma provides a bound on the subsequent function value with each step that holds regardless of the properties of the generated stochastic derivative estimates.

\blemma\label{lem.itr}
  Suppose Assumption~\ref{ass.f} holds.  For all $k \in \N{}$, for any $(g_k,H_k)$, one finds
  \bequationNN
    f(x_{k+1}) \leq f(x_k) + g_k^Ts_k + \half s_k^TH_ks_k + (\nabla f(x_k) - g_k)^Ts_k + \half (L_g + \|H_k\|_2) \|s_k\|_2^2.
  \eequationNN
\elemma
\proof{Proof.}
  Since $f$ is twice continuously differentiable under Assumption~\ref{ass.f}, it follows by Taylor's theorem that there exists $\xhat_k$ on the line segment $[x_k,x_{k+1}]$ such that 
  \bequationNN
    \baligned
      f(x_{k+1}) - f(x_k)
        &= \nabla f(x_k)^Ts_k + \half s_k^T\nabla^2 f(\xhat_k)s_k \\
        &= g_k^Ts_k + \half s_k^TH_ks_k + (\nabla f(x_k) - g_k)^Ts_k + \half s_k^T(\nabla^2 f(\xhat_k)-H_k) s_k.
    \ealigned
  \eequationNN
  Then, since the Cauchy-Schwarz and triangle inequalities together imply with \eqref{eq.Hessian_bounded} that
  \bequationNN
    s_k^T(\nabla^2 f(\xhat_k)-H_k)s_k \leq \|\nabla^2 f(\xhat_k)-H_k\|_2 \|s_k\|_2^2 \leq (L_g + \|H_k\|_2)\|s_k\|_2^2,
  \eequationNN
  the desired result follows. \ifthenelse{\technicalreport=0}{\Halmos}{}
\endproof

Our next lemma is a Cauchy decrease result on the reduction in a quadratic model of the objective function yielded by each computed step.  This type of result is standard in the literature on trust region methods, so we state it without a detailed proof.

\blemma\label{lem.Cauchy}
  For all $k \in \N{}$, for any $(g_k,H_k)$, one finds
  \bequationNN
    g_k^T s_k + \half s_k^T H_k s_k \leq -\half \|g_k\|_2 \min \left\{\Delta_k, \frac{\|g_k\|_2}{\|H_k\|_2} \right\}.
  \eequationNN
\elemma
\proof{Proof.}
  The result follows in the standard manner for Cauchy decrease as in the trust region method literature (see, e.g., \cite[Corollary~6.3.2]{ConnGoulToin00} or \cite[Lemma~4.3]{NoceWrig06}) since, whichever subproblem is solved to compute the step $s_k$, \ref{alg.trish} requires that at least Cauchy decrease is attained. \ifthenelse{\technicalreport=0}{\Halmos}{}
\endproof

We also make use of a second Cauchy decrease result, stated below as our third lemma.  This lemma is useful only when one adds an additional assumption that the norm of the stochastic Hessian estimate is sufficiently small.  We shall add such an assumption for one of our main theorems.  (The proof the lemma follows using a similar argument as in the standard proof for Lemma~\ref{lem.Cauchy}, but with alternative final steps.)

\blemma\label{lem.Cauchy_2}
  For all $k \in \N{}$, for any $(g_k,H_k)$, one finds
  \bequationNN
    g_k^T s_k + \half s_k^T H_k s_k \leq -\min\left\{\Delta_k\|g_k\|_2 - \half \Delta_k^2\|H_k\|_2, \half \frac{\|g_k\|_2^2}{\|H_k\|_2}\right\}.
  \eequationNN
\elemma
\proof{Proof.}
  Using standard analysis for the Cauchy point (see, e.g., \cite[Lemma~4.3]{NoceWrig06}, one has that the Cauchy point lies in the interior of the trust region constraint if $\|g_k\|_2^3 \leq \Delta_k g_k^TH_kg_k$, and lies on the boundary of the trust region constraint otherwise.  If the Cauchy point lies in the interior,  then it is given by $s_k^C := - (\|g_k\|_2^2/g_k^TH_kg_k) g_k$, meaning that, by the Cauchy-Schwarz inequality, the step $s_k$ must satisfy
  \bequationNN
    g_k^Ts_k + \half s_k^TH_ks_k \leq g_k^Ts_k^C + \half {s_k^C}^TH_ks_k^C = -\half \frac{\|g_k\|_2^4}{g_k^TH_kg_k} \leq -\half \frac{\|g_k\|_2^2}{\|H_k\|_2}.
  \eequationNN
  On the other hand, if the Cauchy point lies on the boundary of the trust region constraint, then it is given by $s_k^C := -(\Delta_k/\|g_k\|_2) g_k$ and the step $s_k$ must satisfy
  \bequationNN
    g_k^Ts_k + \half s_k^TH_ks_k \leq g_k^Ts_k^C + \half {s_k^C}^TH_ks_k^C = -\Delta_k \|g_k\|_2 + \half \Delta_k^2 \frac{g_k^TH_kg_k}{\|g_k\|_2^2} \leq -\Delta_k \|g_k\|_2 + \Delta_k^2 \|H_k\|_2.
  \eequationNN
  The result follows by combining the conclusions of these two cases. \ifthenelse{\technicalreport=0}{\Halmos}{}
\endproof

Our next lemma shows that if the stepsize parameter $\alpha_k$ is sufficiently small relative to a quantity involving $\gamma_{1,k}$, $\gamma_{2,k}$, and $\|H_k\|_2$, then the expected reduction in the objective function value with each step is bounded by a function of the expected squared norm of the stochastic gradient estimate, the variance of the stochastic gradient estimate, and the algorithm parameters.  The bound on the reduction proved here will be refined in various ways later in our analysis as we consider the behavior of the algorithm under different sets of assumptions on the derivative estimates and on the stepsize and parameter sequences.

\blemma\label{lem.fundamental}
  Suppose that Assumption~\ref{ass.f} holds and that, for all $k \in \N{}$,
  \bequation\label{eq.alpha}
    0 < \alpha_k \leq \frac{\gamma_{2,k}}{4\gamma_{1,k}^2(L_g + \|H_k\|_2)}.
  \eequation
  Then, for all $k \in \N{}$, one finds
  \bequationNN
    \E_k[f(x_{k+1})] \leq f(x_k) - \frac18 \gamma_{2,k} \alpha_k \E_k[\|g_k\|_2^2] + \frac{\gamma_{1,k}^2}{\gamma_{2,k}} \alpha_k \E_k[\|\nabla f(x_k) - g_k\|_2^2].
  \eequationNN
\elemma
\proof{Proof.}
  We divide the proof according to the three cases defined on page~\pageref{case.1}.
  \bitemize
    \item[\ref{case.1}.] By Lemma~\ref{lem.Cauchy}, it follows in this case that
    \bequationNN
      g_k^Ts_k + \half s_k^TH_ks_k \leq -\half \|g_k\|_2 \min\left\{ \gamma_{1,k}\alpha_k\|g_k\|_2 , \frac{\|g_k\|_2}{\|H_k\|_2} \right\}.
    \eequationNN
    Since \eqref{eq.alpha} ensures $\gamma_{1,k}\alpha_k \leq \tfrac{\gamma_{2,k}}{4\gamma_{1,k}(L_g + \|H_k\|_2)} \leq \tfrac{1}{4(L_g + \|H_k\|_2)} \leq \tfrac{1}{\|H_k\|_2}$, this implies that
    \bequationNN
      g_k^Ts_k + \half s_k^TH_ks_k \leq -\half \gamma_{1,k}\alpha_k \|g_k\|_2^2.
    \eequationNN
    Combining this with the result of Lemma~\ref{lem.itr}, the Cauchy-Schwarz inequality, and the fact that $\|s_k\|_2 \leq \gamma_{1,k}\alpha_k\|g_k\|_2$ in this case, one finds that
    \begin{align}
          &\ f(x_{k+1}) - f(x_k) \nonumber \\
      \leq&\ g_k^Ts_k + \half s_k^TH_ks_k + (\nabla f(x_k) - g_k)^Ts_k + \half (L_g + \|H_k\|_2) \|s_k\|_2^2 \nonumber \\
      \leq&\ -\half \gamma_{1,k}\alpha_k \|g_k\|_2^2 + \|\nabla f(x_k) - g_k\|_2\|s_k\|_2 + \half (L_g + \|H_k\|_2) \|s_k\|_2^2 \nonumber \\
      \leq&\ -\half \gamma_{1,k}\alpha_k \|g_k\|_2^2 + \gamma_{1,k} \alpha_k \|\nabla f(x_k) - g_k\|_2\|g_k\|_2 + \half \gamma_{1,k}^2 \alpha_k^2 (L_g + \|H_k\|_2) \|g_k\|_2^2. \label{eq.for_next}
    \end{align}
    Since $0 \leq (\thalf \|g_k\|_2 - \|\nabla f(x_k) - g_k\|_2)^2 = \tfrac14\|g_k\|_2^2 - \|\nabla f(x_k) - g_k\|_2\|g_k\|_2 + \|\nabla f(x_k) - g_k\|_2^2$ and since \eqref{eq.alpha} implies $\gamma_{1,k}\alpha_k \leq \tfrac{\gamma_{2,k}}{4\gamma_{1,k}(L_g + \|H_k\|_2)} \leq \tfrac{1}{4(L_g + \|H_k\|_2)}$, one finds that
    \bequationNN
      \baligned
            &\ f(x_{k+1}) - f(x_k) \\
        \leq&\ -\half \gamma_{1,k} \alpha_k \|g_k\|_2^2 + \gamma_{1,k} \alpha_k \(\frac14\|g_k\|_2^2 + \|\nabla f(x_k) - g_k\|_2^2\) + \half \gamma_{1,k}^2 \alpha_k^2 (L_g + \|H_k\|_2) \|g_k\|_2^2 \\
        \leq&\ -\frac18 \gamma_{1,k} \alpha_k \|g_k\|_2^2 + \gamma_{1,k} \alpha_k \|\nabla f(x_k) - g_k\|_2^2,
      \ealigned
    \eequationNN
    which implies the desired conclusion since $\gamma_{1,k} \geq \gamma_{2,k}$.
    \item[\ref{case.2}.] By Lemma~\ref{lem.Cauchy} and since in this case one has $\gamma_{2,k} \|g_k\|_2 \leq 1$, it follows that 
    \bequationNN
      g_k^Ts_k + \half s_k^TH_ks_k \leq -\half \|g_k\|_2 \min \left\{\alpha_k, \frac{\|g_k\|_2}{\|H_k\|_2}\right\} \leq -\half \|g_k\|_2 \min \left\{\gamma_{2,k} \alpha_k \|g_k\|, \frac{\|g_k\|_2}{\|H_k\|} \right\}.
    \eequationNN
    Since \eqref{eq.alpha} ensures $\gamma_{2,k}\alpha_k \leq \gamma_{1,k}\alpha_k \leq \tfrac{\gamma_{2,k}}{4\gamma_{1,k}(L_g + \|H_k\|_2)} \leq \tfrac{1}{4(L_g + \|H_k\|_2)} \leq \tfrac{1}{\|H_k\|_2}$, this implies that
    \bequationNN
      g_k^Ts_k + \half s_k^TH_ks_k \leq -\half \gamma_{2,k} \alpha_k \|g_k\|_2^2.
    \eequationNN
    Combining this with the result of Lemma~\ref{lem.itr}, the Cauchy-Schwarz inequality, and the fact that $\|s_k\|_2 \leq \alpha_k$ in this case, one finds that
    \bequationNN
      \baligned
            &\ f(x_{k+1}) - f(x_k) \\
        \leq&\ g_k^Ts_k + \half s_k^TH_ks_k + (\nabla f(x_k) - g_k)^Ts_k + \half (L_g + \|H_k\|_2) \|s_k\|_2^2 \\
        \leq&\ -\half \gamma_{2,k}\alpha_k \|g_k\|_2^2 + \|\nabla f(x_k) - g_k\|_2\|s_k\|_2 + \half (L_g + \|H_k\|_2) \|s_k\|_2^2 \\
        \leq&\ -\half \gamma_{2,k}\alpha_k \|g_k\|_2^2 + \alpha_k \|\nabla f(x_k) - g_k\|_2 + \half \alpha_k^2 (L_g + \|H_k\|_2).
      \ealigned
    \eequationNN
    Since $0 \leq \tfrac{\gamma_{2,k}}{\gamma_{1,k}^2}(\thalf - \tfrac{\gamma_{1,k}^2}{\gamma_{2,k}} \|\nabla f(x_k) - g_k\|_2)^2 = \tfrac{\gamma_{2,k}}{4\gamma_{1,k}^2} - \|\nabla f(x_k) - g_k\|_2 + \tfrac{\gamma_{1,k}^2}{\gamma_{2,k}} \|\nabla f(x_k) - g_k\|_2^2$ and since $1 \leq \gamma_{1,k} \|g_k\|_2$ in this case, the above and \eqref{eq.alpha} imply the desired conclusion that
    \bequationNN
      \baligned
            &\ f(x_{k+1}) - f(x_k) \\
        \leq&\ -\half \gamma_{2,k} \alpha_k \|g_k\|_2^2 + \alpha_k \(\frac{\gamma_{2,k}}{4\gamma_{1,k}^2} + \frac{\gamma_{1,k}^2}{\gamma_{2,k}} \|\nabla f(x_k) - g_k\|_2^2\) + \half \alpha_k^2 (L_g + \|H_k\|_2) \\
           =&\ -\half \gamma_{2,k} \alpha_k \|g_k\|_2^2 + \frac{1}{4}\gamma_{2,k}\alpha_k\|g_k\|_2^2 + \frac{\gamma_{1,k}^2}{\gamma_{2,k}}\alpha_k \|\nabla f(x_k) - g_k\|_2^2 + \half \gamma_{1,k}^2 \alpha_k^2 (L_g + \|H_k\|_2) \|g_k\|_2^2 \\
        \leq&\ -\frac18 \gamma_{2,k} \alpha_k \|g_k\|_2^2 + \frac{\gamma_{1,k}^2}{\gamma_{2,k}} \alpha_k \|\nabla f(x_k) - g_k\|_2^2.
      \ealigned
    \eequationNN
    \item[\ref{case.3}.] The proof follows in the same manner as the proof for \ref{case.1}, using $\gamma_{2,k} \leq \gamma_{1,k}$.
  \eitemize
  The desired conclusion follows by combining the results for the three cases. \ifthenelse{\technicalreport=0}{\Halmos}{}
\endproof

Our final fundamental lemma proves a similar type of bound on the expected reduction in the objective function as in the preceding lemma, except that it can offer a stronger bound when the difference $\gamma_{1,k} - \gamma_{2,k}$ is proportional to $\alpha_k$ and there is an appropriate balance between the stepsize $\alpha_k$ and the norm of the stochastic Hessian estimate.  (Note that to ensure the bound on $\|H_k\|_2$ that is required for the lemma, one might need to scale $H_k$, causing $\E_k[H_k] \neq \nabla^2 f(x_k)$.  This might not seem ideal, but as is known in the deterministic optimization literature, it still allows one to incorporate some (approximate) second-order information, which can be beneficial in practice.)  We consider the behavior of the algorithm in such situations in one of our main theorems.

\blemma\label{lem.fundamental_2}
  Suppose that Assumption~\ref{ass.f} holds and that, for all $k \in \N{}$ and some $\eta \in \R{}_{>0}$,
  \bequation\label{eq.alpha_2}
    \baligned
      0 < \alpha_k &\leq \min\left\{\frac{\gamma_{2,k}}{4\gamma_{1,k}^2(L_g + \|H_k\|_2)}, \frac{1}{6\eta + 2\gamma_{1,k}(L_g + \|H_k\|_2)} \right\}, \\
      \|H_k\|_2 &\leq \frac{\eta}{2\gamma_{1,k}},\ \ \text{and}\ \ \gamma_{1,k} - \gamma_{2,k} = \half \eta \gamma_{1,k} \alpha_k.
    \ealigned
  \eequation
  $($For one thing, this ensures that \eqref{eq.alpha} holds for all $k \in \N{}$.$)$  Then, for all $k \in \N{}$, one finds
  \bequationNN
    \baligned
      \E_k[f(x_{k+1})] \leq f(x_k) &- \frac14 \gamma_{2,k} \alpha_k \|\nabla f(x_k)\|_2^2 \\
      &+ \half \(3\eta + \gamma_{1,k} (L_g + \|H_k\|_2)\) \gamma_{1,k} \alpha_k^2 \E_k[\|\nabla f(x_k) - g_k\|_2^2].
    \ealigned
  \eequationNN
\elemma
\proof{Proof.}
  We divide the proof according to the three cases defined on page~\pageref{case.1}.
  \bitemize
    \item[\ref{case.1}.] By \eqref{eq.alpha}, it follows that $\gamma_{1,k} \alpha_k \leq \frac{\gamma_{2,k}}{4\gamma_{1,k}(L_g + \|H_k\|_2)} \leq \frac{1}{4(L_g + \|H_k\|_2)} \leq \frac{1}{2\|H_k\|_2}$ for all $k \in \N{}$, meaning that for all $k \in \N{}$ one finds in this case that
    \bequationNN
      \Delta_k\|g_k\|_2 - \half \Delta_k^2 \|H_k\|_2 = \gamma_{1,k} \alpha_k \|g_k\|_2^2 - \half \gamma_{1,k}^2 \alpha_k^2 \|g_k\|_2^2 \|H_k\|_2 \leq \half \frac{\|g_k\|_2^2}{\|H_k\|_2},
    \eequationNN
    while at the same time $\half \|H_k\|_2 \leq 2\|H_k\|_2 \leq \frac{\eta}{\gamma_{1,k}}$, meaning for all $k \in \N{}$ that
    \bequationNN
      \baligned
        \Delta_k\|g_k\|_2 - \half \Delta_k^2 \|H_k\|_2 
          &=    \gamma_{1,k} \alpha_k \|g_k\|_2^2 - \half \gamma_{1,k}^2 \alpha_k^2 \|g_k\|_2^2 \|H_k\|_2 \\
          &\geq \gamma_{1,k} \alpha_k \|g_k\|_2^2 - \gamma_{1,k} \alpha_k^2 \eta \|g_k\|_2^2 = (1 - \alpha_k \eta) \gamma_{1,k} \alpha_k \|g_k\|_2^2.
      \ealigned
    \eequationNN
    (Observe that \eqref{eq.alpha_2} ensures that $\alpha < \frac{1}{\eta}$, meaning that $1 - \alpha_k\eta > 0$.) Combining these facts with the results of Lemmas~\ref{lem.itr} and \ref{lem.Cauchy_2}, the Cauchy-Schwarz inequality, and the fact that $\|s_k\|_2 \leq \gamma_{1,k}\alpha_k\|g_k\|_2$ in this case, one finds that
    \bequationNN
      \baligned
            &\ f(x_{k+1}) - f(x_k) \\
        \leq&\ g_k^Ts_k + \half s_k^TH_ks_k + (\nabla f(x_k) - g_k)^Ts_k + \half (L_g + \|H_k\|_2) \|s_k\|_2^2 \\
        \leq&\ - (1 - \alpha_k\eta) \gamma_{1,k} \alpha_k \|g_k\|_2^2 + \|\nabla f(x_k) - g_k\|_2 \|s_k\|_2 + \half (L_g + \|H_k\|_2) \|s_k\|_2^2 \\
        \leq&\ - (1 - \alpha_k\eta) \gamma_{1,k} \alpha_k \|g_k\|_2^2 + \gamma_{1,k} \alpha_k \|\nabla f(x_k) - g_k\|_2 \|g_k\|_2 + \half \gamma_{1,k}^2 \alpha_k^2 (L_g + \|H_k\|_2) \|g_k\|_2^2.
      \ealigned
    \eequationNN
    Since $0 \leq \half(\|g_k\|_2 - \|\nabla f(x_k) - g_k\|_2)^2 = \half \|g_k\|_2^2 - \|\nabla f(x_k) - g_k\|_2 \|g_k\|_2 + \half \|\nabla f(x_k) - g_k\|_2^2$,
    \bequationNN
      \baligned
            &\ f(x_{k+1}) - f(x_k) \\
        \leq&\ -\(1 - \alpha_k \eta - \half \gamma_{1,k} \alpha_k (L_g + \|H_k\|_2)\) \gamma_{1,k} \alpha_k \|g_k\|_2^2 + \half \gamma_{1,k} \alpha_k (\|g_k\|_2^2 + \|\nabla f(x_k) - g_k\|_2^2) \\
           =&\ -\(\half - \alpha_k \eta - \half \gamma_{1,k} \alpha_k (L_g + \|H_k\|_2)\) \gamma_{1,k} \alpha_k \|g_k\|_2^2 + \half \gamma_{1,k} \alpha_k \|\nabla f(x_k) - g_k\|_2^2,
      \ealigned
    \eequationNN
    which along with \eqref{eq.unbiased_consequence} (applied twice) implies that
    \bequationNN
      \baligned
            &\ \E_k[f(x_{k+1})] - f(x_k) \\
        \leq&\ -\(\half - \alpha_k \eta - \half \gamma_{1,k} \alpha_k (L_g + \|H_k\|_2)\) \gamma_{1,k} \alpha_k \E_k[\|g_k\|_2^2] + \half \gamma_{1,k} \alpha_k (-\|\nabla f(x_k)\|_2^2 + \E_k[\|g_k\|_2^2]) \\
           =&\ -\half \gamma_{1,k} \alpha_k \|\nabla f(x_k)\|_2^2 + \(\eta + \half \gamma_{1,k} (L_g + \|H_k\|_2)\) \gamma_{1,k} \alpha_k^2 \E_k[\|g_k\|_2^2] \\
           =&\ -\half \gamma_{1,k} \alpha_k \|\nabla f(x_k)\|_2^2 + \(\eta + \half \gamma_{1,k} (L_g + \|H_k\|_2)\) \gamma_{1,k} \alpha_k^2 (\|\nabla f(x_k)\|_2^2 + \E_k[\|\nabla f(x_k) - g_k\|_2^2]) \\
           =&\ -\(\half - \(\eta + \half \gamma_{1,k} (L_g + \|H_k\|_2)\) \alpha_k\) \gamma_{1,k} \alpha_k \|\nabla f(x_k)\|_2^2 \\
            &\qquad + \(\eta + \half \gamma_{1,k} (L_g + \|H_k\|_2)\) \gamma_{1,k} \alpha_k^2 \E_k[\|\nabla f(x_k) - g_k\|_2^2].
      \ealigned
    \eequationNN
    Hence, with the inequality above, the desired result follows in this case due to the upper bound imposed on $\alpha_k$ and the fact that $\gamma_{1,k} \geq \gamma_{2,k}$ for all $k \in \N{}$.
    \item[\ref{case.2}.] Under the conditions of the lemma, one has $\alpha_k \leq \frac{1}{\eta}$ and $2 \|H_k\|_2 \leq \frac{\eta}{\gamma_{1,k}}$.  In addition, in this case, one has $\gamma_{1,k} \|g_k\|_2 \geq 1$.  These facts combined imply that
    \bequationNN
      \baligned
        \Delta_k\|g_k\|_2 - \half \Delta_k^2 \|H_k\|_2 &= \alpha_k \|g_k\|_2 - \half \alpha_k^2 \|H_k\|_2 \leq \frac{\|g_k\|_2}{\eta} \\
        \text{while}\ \ \half \frac{\|g_k\|_2^2}{\|H_k\|_2} &\geq \frac{\gamma_{1,k} \|g_k\|_2^2}{\eta} \geq \frac{\|g_k\|_2}{\eta}. \\
      \ealigned
    \eequationNN
    By Lemma~\ref{lem.Cauchy_2} and the facts that $\half \|H_k\|_2 \leq 2 \|H_k\|_2 \leq \frac{\eta}{\gamma_{1,k}}$ and $\gamma_{1,k}\|g_k\|_2 \geq 1$, it follows that
    \bequationNN
      \baligned
        g_k^Ts_k + \half s_k^TH_ks_k
          &\leq -\Delta_k \|g_k\|_2 + \half \Delta_k^2 \|H_k\|_2 \\
          &=    -\alpha_k \|g_k\|_2 + \half \alpha_k^2 \|H_k\|_2 \\
          &\leq -(1 - \alpha_k\eta) \alpha_k \|g_k\|_2.
      \ealigned
    \eequationNN
    Combining this fact with the results of Lemmas~\ref{lem.itr} and \ref{lem.Cauchy_2}, the Cauchy-Schwarz inequality, and the facts that $\gamma_{2,k}\|g_k\|_2 \leq 1$, $\gamma_{1,k}\|g_k\|_2 \geq 1$, and $\|s_k\|_2 \leq \alpha_k$ in this case, one finds that
    \bequationNN
      \baligned
            &\ f(x_{k+1}) - f(x_k) \\
        \leq&\ g_k^Ts_k + \half s_k^TH_ks_k + (\nabla f(x_k) - g_k)^Ts_k + \half (L_g + \|H_k\|_2) \|s_k\|_2^2 \\
        \leq&\ -(1 - \alpha_k\eta)\alpha_k\|g_k\|_2 + \|\nabla f(x_k) - g_k\|_2 \|s_k\|_2 + \half(L_g + \|H_k\|_2) \|s_k\|_2^2 \\
        \leq&\ -(1 - \alpha_k\eta)\gamma_{2,k}\alpha_k\|g_k\|_2^2 + \alpha_k \|\nabla f(x_k) - g_k\|_2 + \half \alpha_k^2 (L_g + \|H_k\|_2) \\
        \leq&\ -(1 - \alpha_k\eta)\gamma_{2,k}\alpha_k\|g_k\|_2^2 + \gamma_{1,k} \alpha_k \|\nabla f(x_k) - g_k\|_2 \|g_k\|_2 + \half \gamma_{1,k}^2 \alpha_k^2 (L_g + \|H_k\|_2) \|g_k\|_2^2.
      \ealigned
    \eequationNN
    Since $0 \leq \half(\|g_k\|_2 - \|\nabla f(x_k) - g_k\|_2)^2 = \half \|g_k\|_2^2 - \|\nabla f(x_k) - g_k\|_2 \|g_k\|_2 + \half \|\nabla f(x_k) - g_k\|_2^2$,
    \bequationNN
      \baligned
            &\ f(x_{k+1}) - f(x_k) \\
        \leq&\ -\((1 - \alpha_k\eta)\gamma_{2,k} - \half \gamma_{1,k}^2 \alpha_k(L_g + \|H_k\|_2)\) \alpha_k \|g_k\|_2^2 + \half \gamma_{1,k} \alpha_k \(\|g_k\|_2^2 + \|\nabla f(x_k) - g_k\|_2^2\) \\
           =&\ -\((1 - \alpha_k\eta)\gamma_{2,k} - \half \gamma_{1,k} - \half \gamma_{1,k}^2 \alpha_k(L_g + \|H_k\|_2)\) \alpha_k \|g_k\|_2^2 + \half \gamma_{1,k} \alpha_k \|\nabla f(x_k) - g\|_2^2,
      \ealigned
    \eequationNN
    which along with \eqref{eq.unbiased_consequence} (applied twice) implies that
    \bequationNN
      \baligned
            &\ \E_k[f(x_{k+1})] - f(x_k) \\
        \leq&\ -\((1 - \alpha_k\eta)\gamma_{2,k} - \half \gamma_{1,k} - \half \gamma_{1,k}^2 \alpha_k(L_g + \|H_k\|_2)\) \alpha_k \E_k[\|g_k\|_2^2] \\
            &\qquad + \half \gamma_{1,k} \alpha_k (-\|\nabla f(x_k)\|_2^2 + \E_k[\|g_k\|_2^2]) \\
           =&\ -\half \gamma_{1,k} \alpha_k \|\nabla f(x_k)\|_2^2 + \(\gamma_{1,k} - \gamma_{2,k} + \(\eta\gamma_{2,k} + \half \gamma_{1,k}^2 (L_g + \|H_k\|_2)\)\alpha_k\) \alpha_k \E_k[\|g_k\|_2^2] \\
        \leq&\ -\half \gamma_{1,k} \alpha_k \|\nabla f(x_k)\|_2^2 \\
            &\qquad + \(\gamma_{1,k} - \gamma_{2,k} + \(\eta\gamma_{2,k} + \half \gamma_{1,k}^2 (L_g + \|H_k\|_2)\)\alpha_k\) \alpha_k (\|\nabla f(x_k)\|_2^2 + \E_k[\|\nabla f(x_k) - g_k\|_2^2]) \\
           =&\ -\(\half \gamma_{1,k} - (\gamma_{1,k} - \gamma_{2,k}) - \(\eta\gamma_{2,k} + \half \gamma_{1,k}^2 (L_g + \|H_k\|_2)\)\alpha_k\) \alpha_k \|\nabla f(x_k)\|_2^2 \\
            &\qquad + \(\gamma_{1,k} - \gamma_{2,k} + \(\eta\gamma_{2,k} + \half \gamma_{1,k}^2 (L_g + \|H_k\|_2)\)\alpha_k\) \alpha_k \E_k[\|\nabla f(x_k) - g_k\|_2^2] \\
           =&\ -\(\half - \half \(3\eta - \eta^2 \alpha_k + \gamma_{1,k} (L_g + \|H_k\|_2)\)\alpha_k\) \gamma_{1,k} \alpha_k \|\nabla f(x_k)\|_2^2 \\
            &\qquad + \half\(3\eta - \eta^2 \alpha_k + \gamma_{1,k} (L_g + \|H_k\|_2)\) \gamma_{1,k} \alpha_k^2 \E_k[\|\nabla f(x_k) - g_k\|_2^2].
      \ealigned
    \eequationNN
    Hence, the desired result follows for this case, again due to the upper bound on $\alpha_k$ and the fact that $\gamma_{1,k} \geq \gamma_{2,k}$ for all $k \in \N{}$.
    \item[\ref{case.3}.]  The proof for this case follows in the same manner as the proof for \ref{case.1}, where the result for this case has a similar form except with $\gamma_{1,k}$ replaced by $\gamma_{2,k}$.  For the proof, it should be noted that $\gamma_{2,k}\alpha_k \leq \gamma_{1,k}\alpha_k \leq \frac{1}{2\|H_k\|_2}$, $\half \|H_k\|_2 \leq \frac{\eta}{\gamma_{1,k}} \leq \frac{\eta}{\gamma_{2,k}}$, and $\|s_k\|_2 \leq \gamma_{2,k}\alpha_k\|g_k\|_2$.
  \eitemize
  The desired conclusion follows by combining the results for the three cases. \ifthenelse{\technicalreport=0}{\Halmos}{}
\endproof

Now that these fundamental lemmas have been established, which focus on the behavior of the algorithm over a single iteration, we turn to analyzing the behavior of the algorithm over the entire sequence of iterations.  We break our analysis into parts based on different assumptions about the problem function and the stochastic derivative estimates.  For simplicity in much of our analysis, we consider the behavior of the algorithm when the parameter sequences $\{\gamma_{1,k}\}$ and $\{\gamma_{2,k}\}$ are constant.  In such cases, one could prove similar results that allow the sequences not to be constant, as long as they remain within bounded intervals.  We also prove one result showing that, in practice, one might define these sequences to have the same limit point, which in a sense makes the algorithm behave asymptotically like a more standard stochastic Newton-type method.

\subsection{General (nonconvex) objective functions}\label{sec.general}

First, we consider the case when the algorithm is employed to minimize an objective function satisfying only Assumptions~\ref{ass.f} and \ref{ass.gH}, and when the following loose assumption holds about the algorithm parameters, stochastic gradients, and stochastic Hessians.

\bassumption\label{ass.nonconvex.fixed}
  The variance of the stochastic gradient estimates and the sequence of stochastic Hessian estimates are both uniformly bounded in the sense that there exist constants $(M_g,M_H) \in \R{}_{>0} \times \R{}_{>0}$ such that, for all $k \in \N{}$,
  \bequationNN
    \E_k[\|\nabla f(x_k) - g_k\|_2^2] \leq M_g\ \ \text{and}\ \ \|H_k\|_2 \leq M_H.
  \eequationNN
  In addition, $(\alpha_k,\gamma_{1,k},\gamma_{2,k}) = (\alpha,\gamma_1,\gamma_2)$ for all $k \in \N{}$, where $\gamma_1 \geq \gamma_2 > 0$ and
  \bequationNN
    0 < \alpha \leq \frac{\gamma_2}{4\gamma_1^2(L_g + M_H)},
  \eequationNN
  which, in particular, implies that \eqref{eq.alpha} holds for all $k \in \N{}$.
\eassumption

Combining the result of Lemma~\ref{lem.fundamental} with Assumption~\ref{ass.nonconvex.fixed} leads to the following result showing that the expected average squared norm of the gradient at the iterates is bounded.

\btheorem\label{th.nonconvex.fixed}
  Under Assumptions \ref{ass.f}, \ref{ass.gH}, and \ref{ass.nonconvex.fixed}, \ref{alg.trish} yields
  \bsubequations
    \begin{align}
      \E\left[\sum_{k=1}^K \|\nabla f(x_k)\|_2^2\right] \leq&\ \(\frac{8}{\gamma_2 \alpha}\) (f(x_1) - f_{\rm inf}) + K\(\frac{8\gamma_1}{\gamma_2} - 1\)M_g \label{eq.nonconvex.fixed.1} \\
      \text{and}\ \ \E\left[\frac{1}{K} \sum_{k=1}^K \|\nabla f(x_k)\|_2^2\right] \leq&\ \frac{1}{K} \(\frac{8}{\gamma_2 \alpha}\) (f(x_1) - f_{\rm inf}) + \(\frac{8\gamma_1^2}{\gamma_2^2} - 1\)M_g \nonumber \\
      \xrightarrow{K\to\infty} & \(\frac{8\gamma_1^2}{\gamma_2^2} - 1\)M_g. \label{eq.nonconvex.fixed.2}
    \end{align}
  \esubequations
\etheorem
\proof{Proof.}
  Since Assumption~\ref{ass.nonconvex.fixed} ensures that \eqref{eq.alpha} holds for all $k \in \N{}$, it follows that the result of Lemma~\ref{lem.fundamental} holds; hence, with parameters as in Assumption~\ref{ass.nonconvex.fixed}, for all $k \in \N{}$ one has
  \bequationNN
    \E_k [f(x_{k+1})] \leq f(x_k) - \frac18 \gamma_2 \alpha \E_k [\|g_k\|_2^2] + \frac{\gamma_1^2}{\gamma_2} \alpha \E_k [\|\nabla f(x_k) - g_k\|_2^2].
  \eequationNN
  Hence, due to Assumption~\ref{ass.nonconvex.fixed} and \eqref{eq.unbiased_consequence}, it follows for all $k \in \N{}$ that
  \begin{align}
        &\ \E_k [f(x_{k+1})] - f(x_k) \nonumber \\
    \leq&\ - \frac18 \gamma_2 \alpha (\|\nabla f(x_k)\|_2^2 + \E_k[\|\nabla f(x_k) - g_k\|_2^2]) + \frac{\gamma_1^2}{\gamma_2} \alpha \E_k [\|\nabla f(x_k) - g_k\|_2^2] \nonumber \\
       =&\ - \frac18 \gamma_2 \alpha \|\nabla f(x_k)\|_2^2 + \alpha\(\frac{\gamma_1^2}{\gamma_2} - \frac18\gamma_2\)\E_k[\|\nabla f(x_k) - g_k\|_2^2] \nonumber \\
    \leq&\ - \frac18 \gamma_2 \alpha \|\nabla f(x_k)\|_2^2 + \alpha\(\frac{\gamma_1^2}{\gamma_2} - \frac18\gamma_2\)M_g. \label{eq.for_later}
  \end{align}
  Taking total expectation, it follows for all $k \in \N{}$ that
  \bequationNN
    \E[f(x_{k+1})] - \E[f(x_k)] \leq - \frac18 \gamma_2 \alpha \E[\|\nabla f(x_k)\|_2^2] + \alpha\(\frac{\gamma_1^2}{\gamma_2} - \frac18\gamma_2\)M_g,
  \eequationNN
  which implies
  \bequationNN
    \E[\|\nabla f(x_k)\|_2^2] \leq \(\frac{8}{\gamma_2 \alpha}\)(\E[f(x_k)] - \E[f(x_{k+1})]) + \(\frac{8\gamma_1^2}{\gamma_2^2} - 1\)M_g.
  \eequationNN
  Summing this inequality over all $k \in \{1,\dots,K\}$ and using the fact that $f$ is bounded below by $\finf$ yields \eqref{eq.nonconvex.fixed.1}, which, in turn, implies \eqref{eq.nonconvex.fixed.2}. \ifthenelse{\technicalreport=0}{\Halmos}{}
\endproof

Next, we consider the behavior of \ref{alg.trish} when Assumptions~\ref{ass.f} and \ref{ass.gH} hold and when a run of the algorithm satisfies the following assumption involving diminishing stepsizes.

\bassumption\label{ass.nonconvex.diminishing}
  The variance of each stochastic gradient estimate is proportional to the stepsize and the sequence of stochastic Hessian estimates is uniformly bounded in the sense that there exist constants $(M_g,M_H) \in \R{}_{>0} \times \R{}_{>0}$ such that, for all $k \in \N{}$,
  \bequation\label{eq.nonconvex.diminishing.variance}
    \E_k[\|\nabla f(x_k) - g_k\|_2^2] \leq M_g \alpha_k\ \ \text{and}\ \ \|H_k\|_2 \leq M_H.
  \eequation
  In addition, $(\gamma_{1,k},\gamma_{2,k}) = (\gamma_1,\gamma_2)$ for all $k \in \N{}$ where $\gamma_1 \geq \gamma_2 > 0$, and
  \bequationNN
    \{\alpha_k\} = \left\{\frac{a}{b+k}\right\}\ \ \text{for some $(a,b) \in \R{}_{>0} \times \R{}_{>0}$}
  \eequationNN
  such that \eqref{eq.alpha} holds for all $k \in \N{}$.
\eassumption

Under this assumption, which is stronger than Assumption~\ref{ass.nonconvex.fixed}, we obtain the following result, which, not surprisingly, is stronger than the result in Theorem~\ref{th.nonconvex.fixed}.

\btheorem\label{th.nonconvex.diminishing}
  Under Assumptions \ref{ass.f}, \ref{ass.gH}, and \ref{ass.nonconvex.diminishing}, \ref{alg.trish} yields
  \bsubequations
    \begin{align}
      \lim_{K\to\infty} \E\left[ \sum_{k=1}^K \alpha_k \|\nabla f(x_k)\|^2 \right] &< \infty \label{eq.nonconvex.diminishing.1} \\
      \text{and}\ \ \E\left[ \frac{1}{\sum_{k=1}^K \alpha_k} \sum_{k=1}^K \alpha_k \|\nabla f(x_k)\|^2 \right] &\xrightarrow{K\to\infty} 0. \label{eq.nonconvex.diminishing.2}
    \end{align}
  \esubequations
  In addition, it follows that
  \bequation\label{eq.nonconvex.diminishing.3}
    \liminf_{k \rightarrow \infty} \|f(x_k)\|_2^2 = 0\ \ \text{with probability 1}.
  \eequation
\etheorem
\proof{Proof.}
  Following the same arguments as in the proof of Theorem~\ref{th.nonconvex.fixed}, for all $k \in \N{}$,
  \bequation\label{eq.nonconvex.diminishing.reduction}
    \E_k[f(x_{k+1})] \leq f(x_k) - \frac18 \gamma_2 \alpha_k \|\nabla f(x_k)\|_2^2 + \(\frac{\gamma_1^2}{\gamma_2} - \frac18\gamma_2\) M_g \alpha_k^2,
  \eequation
  which, taking total expectation, implies for all $k \in \N{}$ that
  \bequationNN
    \E[f(x_{k+1})] - \E[f(x_k)] \leq -\frac18 \gamma_2 \alpha_k \E[\|\nabla f(x_k)\|_2^2] + \(\frac{\gamma_1^2}{\gamma_2} - \frac18\gamma_2\) M_g \alpha_k^2.
  \eequationNN
  Rearranging terms and summing over all $k \in \{1,\dots,K\}$, it follows that
  \bequation\label{eq.pre}
    \frac{1}{8} \gamma_2 \sum_{k=1}^K \alpha_k \E[\|\nabla f(x_k)\|_2^2] \leq \sum_{k=1}^K \(\E[f(x_k)] - \E[f(x_{k+1})]\) + \(\frac{\gamma_1^2}{\gamma_2} - \frac18\gamma_2\) M_g \sum_{k=1}^K \alpha_k^2.
  \eequation
  Since $\sum_{k=1}^K \(\E[f(x_k)] - \E[f(x_{k+1})]\) \leq f(x_1) - \finf < \infty$ for any $K \in \N{}$ and since Assumption~\ref{ass.nonconvex.diminishing} implies $\sum_{k=1}^\infty \alpha_k^2 < \infty$, it follows from \eqref{eq.pre} that \eqref{eq.nonconvex.diminishing.1} holds.  Moreover, dividing \eqref{eq.pre} by $\sum_{k=1}^K \alpha_k$ and since Assumption~\ref{ass.nonconvex.diminishing} implies $\sum_{k=1}^\infty \alpha_k = \infty$, it follows that \eqref{eq.nonconvex.diminishing.2} holds.
  
  Let us now prove \eqref{eq.nonconvex.diminishing.3}.  Defining the scalars $\beta_1 := \tfrac18 \gamma_2$ and $\beta_2 := \(\frac{\gamma_1^2}{\gamma_2} - \frac18\gamma_2\) M_g$, it follows from \eqref{eq.nonconvex.diminishing.reduction} that, for all $k \in \N{}$, the expected reduction in $f$ satisfies
  \bequationNN
    \baligned
      \E_k[f(x_{k+1})] &\leq f(x_k) - \beta_1 \alpha_k \|\nabla f(x_k)\|^2 + \beta_2 \alpha_k^2 \\
      \implies\ \ \E_k[f(x_{k+1})] + \beta_2 \sum_{i=k+1}^{\infty} \alpha_i^2 &\leq f(x_k) - \beta_1 \alpha_k \|\nabla f(x_k)\|^2 + \beta_2 \sum_{i=k}^{\infty} \alpha_i^2.
    \ealigned
  \eequationNN
  Considering the stochastic processes $\{p_k\}$ and $\{q_k\}$, where, for all $k \in \N{}$,
  \bequationNN
    p_k := \beta_1\alpha_k\|\nabla f(x_k)\|_2^2\ \ \text{and}\ \ q_k := f(x_k) + \beta_2\sum_{i=k}^{\infty} \alpha_i^2,
  \eequationNN
  it follows from above that, for all $k \in \N{}$,
  \bequation\label{eq.pq}
    \E_k [q_{k+1}-\finf] \leq q_k - \finf - p_k.
  \eequation
  One finds from this relationship that $\E[q_k-\finf] < \infty$ and $\E_k[q_{k+1}-\finf] \leq q_k - \finf$ for all $k \in \N{}$, which with $q_k-\finf\geq0$ for all $k \in \N{}$ implies that $\{q_k - \finf\}$ is a nonnegative supermartingale.  This implies (see, e.g., \cite{Durr19} and similar use in \cite{WangMaGoldLiu17}) that there exists $q$ such that $\lim_{k \rightarrow \infty} q_k = q$ with probability 1 and $\E[q] \leq \E[q_1]$.  From \eqref{eq.pq}, one finds that $\E[p_k] \leq \E[q_k] - \E[q_{k+1}]$, from which it follows that
  \bequation\label{eq.p_bound}
    \E\left[\sum_{k=1}^\infty p_k\right] < \infty\ \ \implies\ \ \sum_{k=1}^\infty \beta_1\alpha_k\|\nabla f(x_k)\|_2^2 = \sum_{k=1}^\infty p_k < \infty\ \ \text{with probability 1}.
  \eequation
  Since $\sum_{k=1}^\infty \alpha_k = \infty$ under Assumption~\ref{ass.nonconvex.diminishing}, the above implies \eqref{eq.nonconvex.diminishing.3}. \ifthenelse{\technicalreport=0}{\Halmos}{}
\endproof

To conclude this section, let us prove a result that in part considers the behavior of the algorithm under the following assumption.

\bassumption\label{ass.nonconvex.diminishing.2}
  The second moment of the stochastic gradient estimates is uniformly bounded in the sense that there exists a constant $M_{g,2} \in \R{}_{>0}$ such that, for all $k \in \N{}$,
  \bequationNN
    \E_k[\|g_k\|_2^2] \leq M_{g,2}.
  \eequationNN
\eassumption

It should be said that Assumption~\ref{ass.nonconvex.diminishing.2} is strong since it implies that the variance of the stochastic gradient estimates is smaller at points at which $\|\nabla f(x_k)\|_2$ is large.  In particular, under Assumptions~\ref{ass.gH} and Assumption~\ref{ass.nonconvex.diminishing.2}, it follows $($recall \eqref{eq.unbiased_consequence}$)$ that
\bequationNN
  \E_k[\|g_k\|_2^2] \leq M_g\ \ \implies\ \ \E_k[\|\nabla f(x_k) - g_k\|_2^2] \leq M_g - \|\nabla f(x_k)\|_2^2.
\eequationNN
That said, if the iterates of the algorithm happen to remain in a region over which $\|\nabla f(\cdot)\|_2$ is bounded, then it is interesting to note that Assumption~\ref{ass.nonconvex.diminishing.2}  leads to the following strong result about the behavior of the algorithm.  (A result similar to the following was proved for a stochastic quasi-Newton method as \cite[Th.~2.6]{WangMaGoldLiu17}, and our proof borrows from that one.  That said, our proof corrects an oversight made in the proof of \cite[Th.~2.6]{WangMaGoldLiu17} when one considers the negation of a statement of the form \eqref{eq.nonconvex.diminishing.4}; in particular, in the negation, one should only assume that the limit does not hold with some positive probability, not with complete certainty.)

\btheorem \label{thm.twocase.dim}
  Under Assumptions \ref{ass.f}, \ref{ass.gH}, \ref{ass.nonconvex.diminishing}, and \ref{ass.nonconvex.diminishing.2}, \ref{alg.trish} yields
  \bequation\label{eq.nonconvex.diminishing.4}
    \lim_{k \rightarrow \infty} \|\nabla f(x_k)\|_2 = 0\ \ \text{with probability 1.}
  \eequation
\etheorem
\proof{Proof.}
  To derive a contradiction, suppose that \eqref{eq.nonconvex.diminishing.4} does not hold, meaning that with some nonzero probability there exists $\epsilon \in (0,\infty)$ and an infinite index set $\Kcal_1 \subseteq \N{}$ such that $\|\nabla f(x_k)\|_2 > \epsilon$ for all $k \in \Kcal_1$.  On the other hand, from Theorem~\ref{th.nonconvex.diminishing} it follows that~\eqref{eq.nonconvex.diminishing.3} holds, meaning that with probability one there exists an infinite index set $\Kcal_2$ such that $\|\nabla f(x_k)\|_2 \leq \thalf \epsilon$ for all $k \in \Kcal_2$.  Together, these facts imply with nonzero probability the existence of index sets $\{m_i\}_{i=1}^\infty \subset \N{}$ and $\{n_i\}_{i=1}^\infty \subset \N{}$ with $m_i < n_i$ for all $i \in \N{}$ such that
  \bequation\label{eq.epsilon}
    \baligned
      \|\nabla f(x_{m_i})\| \geq \epsilon,\ \ \|\nabla f(x_{n_i})\| &< \thalf\epsilon, \\
      \text{and}\ \ \|\nabla f(x_{k})\| &\geq \thalf\epsilon\ \ \text{for all}\ \ k \in \{m_i+1,\cdots,n_i-1\}.
    \ealigned
  \eequation
  
  For the remainder of the proof, let us condition on the event that \eqref{eq.epsilon} holds.  With \eqref{eq.p_bound},
  \bequationNN
    \infty > \sum_{k=1}^{\infty} \alpha_k \|\nabla f(x_k)\|_2^2 \geq \sum_{i=1}^{\infty} \sum_{k=m_i}^{n_i-1} \alpha_k \|\nabla f(x_k)\|_2^2 \geq \epsilon^2 \sum_{i=1}^{\infty} \sum_{k=m_i}^{n_i-1} \alpha_k\ \ \text{with probability 1},
  \eequationNN
  meaning that
  \bequation\label{eq.nonzero}
    \lim_{i\to\infty} \sum_{k=m_i}^{n_i-1} \alpha_k < \infty\ \ \text{with probability 1}.
  \eequation
  Now notice that, for any $k \in \N{}$, for any $(g_k,H_k)$, Assumption~\ref{ass.nonconvex.diminishing.2} implies
  \bequationNN
    \E_k[\|x_{k+1}-x_k\|_2] = \E_k[\|s_k\|] \leq \alpha_k \max\{1,\gamma_1\E_k[\|g_k\|_2]\} = \alpha_k \max\left\{1,\gamma_1\sqrt{M_{g,2}}\right\},
  \eequationNN
  from which it follows that
  \bequationNN
    \E_k[\|x_{n_i}-x_{m_i}\|_2] \leq \max\left\{1,\gamma_1\sqrt{M_{g,2}}\right\} \sum_{k=m_i}^{n_i-1} \alpha_k.
  \eequationNN
  Therefore, with \eqref{eq.nonzero}, one finds that $\lim_{i\to\infty} \|x_{n_i}-x_{m_i}\|_2 = 0$ with probability 1, which with Lipschitz continuity of $\nabla f$ under Assumption~\ref{ass.f} implies that $\lim_{i\to\infty} \|\nabla f(x_{n_i}) - \nabla f(x_{m_i})\|_2 = 0$ with probability 1.  However, this contradictions \eqref{eq.epsilon}. \ifthenelse{\technicalreport=0}{\Halmos}{}
\endproof

\subsection{Objective functions satisfying the Polyak-\L{}ojasiewicz condition}

We now consider when the algorithm is employed to minimize an objective function satisfying Assumptions~\ref{ass.f} and \ref{ass.gH} along with the Polyak-\L{}ojasiewicz (PL) condition.  We state this condition in the form of the following assumption.

\bassumption\label{ass.PL}
  There exists a constant $c \in (0,\infty)$ such that, for all $x \in \R{n}$, one has
  \bequation\label{eq.PL}
    2c(f(x) - \finf) \leq \|\nabla f(x)\|_2^2\ \ \text{for all}\ \ x \in \R{n}.
  \eequation
\eassumption

Functions satisfying Assumption~\ref{ass.PL} include $c$-strongly convex functions, but also other nonconvex functions.  Assumptions~\ref{ass.f} and \ref{ass.PL} combined do not guarantee that $f$ has a minimizer, although they do guarantee that if a stationary point exists then it is a global minimizer with objective value $\finf$.  The PL condition is known as a relatively weak condition under which certain algorithms, such as gradient descent, can enjoy a linear rate of convergence.  In this section, we show that the theoretical properties for \ref{alg.trish} are stronger under the PL condition than they are in the more general situations considered in \S\ref{sec.general}.

Our first result shows that if the variance of the stochastic gradient estimates and the stochastic Hessian estimates are both uniformly bounded and the algorithm is run with certain fixed parameter settings, then the expected optimality gap is bounded above by a sequence that converges linearly to a constant proportional to $M_g/c$.  This result is comparable to one that can be proved for SG with a fixed stepsize, for which the limiting constant is also $\Ocal(M_g/c)$; see \cite[Theorem 4.6]{BottCurtNoce18}.

\btheorem\label{th.PL.fixed}
  Under Assumptions~\ref{ass.f}, \ref{ass.gH}, \ref{ass.nonconvex.fixed}, and \ref{ass.PL}, if $\alpha \leq 4/(\gamma_2c)$, then with
  \bequation\label{eq.theta}
    \theta := 4\(\frac{\gamma_1^2}{\gamma_2^2} - \frac18\)\frac{M_g}{c}
  \eequation  
  \ref{alg.trish} yields
  \bequationNN
    \E[f(x_{K+1})] - \finf \leq \theta + \(1 - \frac14 \gamma_2 c \alpha\)^K (f(x_1) - \finf - \theta) \xrightarrow{K\to\infty} \theta.
  \eequationNN
\etheorem
\proof{Proof.}
  As in the proof of Theorem~\ref{th.nonconvex.fixed} (see \eqref{eq.for_later}), it follows for all $k \in \N{}$ that
  \bequationNN
    \E_k[f(x_{k+1})] \leq f(x_k) - \frac18 \gamma_2 \alpha \|\nabla f(x_k)\|_2^2 + \alpha \(\frac{\gamma_1^2}{\gamma_2} - \frac18\gamma_2\) M_g.
  \eequationNN
  Hence, by Assumption~\ref{ass.PL}, it follows for all $k \in \N{}$ that
  \bequationNN
    \E_k[f(x_{k+1})] \leq f(x_k) - \frac14 \gamma_2 c \alpha (f(x_k) - \finf) + \alpha \(\frac{\gamma_1^2}{\gamma_2} - \frac18\gamma_2\) M_g.
  \eequationNN
  Subtracting $\finf$ from both sides and taking total expectation, it follows for all $k \in \N{}$ that
  \bequationNN
    \E[f(x_{k+1})] - \finf \leq \(1 - \frac14 \gamma_2 c \alpha\) (\E[f(x_k)] - \finf) + \alpha \(\frac{\gamma_1^2}{\gamma_2} - \frac18\gamma_2\) M_g.
  \eequationNN
  Therefore, with $\theta$ defined in \eqref{eq.theta}, it follows for all $k \in \N{}$ that
  \bequationNN
    \E[f(x_{k+1})] - \finf - \theta \leq \(1 - \frac14 \gamma_2 c \alpha\) (\E[f(x_k)] - \finf - \theta).
  \eequationNN
  Applying this bound repeatedly for $k \in \{1,\dots,K\}$ yields the desired result. \ifthenelse{\technicalreport=0}{\Halmos}{}
\endproof

Let us now prove, under similar assumptions as in the previous theorem (in particular with respect to the stochastic gradient and Hessian estimates), that \ref{alg.trish} can offer sublinear decrease of the expected optimality gap to zero if the stepsizes vanish along with the differences $\{\gamma_{1,k} - \gamma_{2,k}\}$.  This is the only theorem that we prove in which we consider a case in which $\{\gamma_{1,k}\}$ and $\{\gamma_{2,k}\}$ are not both constant; in particular, we assume $\{\gamma_{1,k}\}$ is constant, but that $\{\gamma_{2,k}\}$ is not.  Other similar results can be proved, say with $\{\gamma_{1,k}\}$ converging to a constant sequence $\{\gamma_{2,k}\}$, or with $\{\gamma_{1,k}\}$ and $\{\gamma_{2,k}\}$ both not constant as long as the sequences remain within a positive interval and the difference sequence is proportional to the stepsize sequence in that $\{\gamma_{1,k} - \gamma_{2,k}\} = \Ocal(\alpha_k)$.

For this theorem only, we consider the following assumption.
\bassumption\label{ass.merging}
  The variance of the stochastic gradient estimates and the sequence of stochastic Hessian estimates are both uniformly bounded in the sense that there exist constants $(M_g,M_H) \in \R{}_{>0} \times \R{}_{>0}$ such that, for all $k \in \N{}$,
  \bequationNN
    \E_k[\|\nabla f(x_k) - g_k\|_2^2] \leq M_g\ \ \text{and}\ \ \|H_k\|_2 \leq M_H.
  \eequationNN
  In addition, $\gamma_{1,k} = \gamma_1 > 0$ for all $k \in \N{}$, and
  \bequationNN
    \{\alpha_k\} = \left\{\frac{a}{b+k}\right\}\ \ \text{and}\ \ \{\gamma_{2,k}\} = \left\{\gamma_1\(1 - \half\eta\alpha_k\)\right\}\ \ \text{for some $(a,b,\eta) \in \R{}_{>0} \times \R{}_{>0} \times \R{}_{>0}$}
  \eequationNN
  such that \eqref{eq.alpha_2} holds for all $k \in \N{}$.
\eassumption

Under this assumption, we prove sublinear decrease of the expected optimality gap.

\btheorem\label{th.PL.sublinear_1}
  Under Assumptions~\ref{ass.f}, \ref{ass.gH}, \ref{ass.PL}, and \ref{ass.merging}, if the pair $(a,b) \in \R{}_{>0} \times \R{}_{>0}$ is chosen such that $\alpha_k \leq \frac{2}{\gamma_{2,1}c}$ for all $k \in \N{}$, then for all $k \in \N{}$ the expected optimality gap satisfies
  \bequation\label{eq.sublinear_1}
    \E[f(x_k)] - \finf \leq \frac{\phi}{b + k},
  \eequation
  where
  \bequation\label{eq.phi}
    \phi := \max\left\{(b+1) (f(x_1)-\finf), \frac{\delta_2 a^2}{\delta_1a-1}\right\} \in (0,\infty),
  \eequation
  with
  \bequation\label{eq.delta}
    \delta_1 := \frac12\gamma_{2,1}c \in \(0,\frac{1}{\alpha}\right]\ \ \text{and}\ \ \delta_2 := \half (3\eta + \gamma_1(L_g + M_H)) \gamma_1 M_g \in (0,\infty).
  \eequation
\etheorem
\proof{Proof.}
  By Lemma~\ref{lem.fundamental_2}, it follows for all $k \in \N{}$ that
  \bequation\label{eq.same_as_next}
    \E_k[f(x_{k+1})] \leq f(x_k) - \frac14 \gamma_{2,1} \alpha_k \|\nabla f(x_k)\|_2^2 + \half (3\eta + \gamma_1(L_g + M_H)) \gamma_1 M_g \alpha_k^2.
  \eequation
    Hence, by Assumption~\ref{ass.PL}, it follows for all $k \in \N{}$ that
  \bequationNN
    \E_k[f(x_{k+1})] \leq f(x_k) - \half \gamma_{2,1} c \alpha_k (f(x_k) - \finf) + \half (3\eta + \gamma_1(L_g + M_H)) \gamma_1 M_g \alpha_k^2.
  \eequationNN
  Subtracting $\finf$ from both sides and taking total expectation, it follows for all $k \in \N{}$ that
  \bequationNN
    \E[f(x_{k+1})] - \finf \leq \(1 - \half \gamma_{2,1} c \alpha_k\) (\E[f(x_k)] - \finf) + \half (3\eta + \gamma_1(L_g + M_H)) \gamma_1 M_g \alpha_k^2.
  \eequationNN
  Let us now prove \eqref{eq.sublinear_1} by induction.  First, for $k=1$, the inequality holds by the definition of~$\phi$ in \eqref{eq.phi}.  Now suppose that \eqref{eq.sublinear_1} holds up to $k \in \N{}$.  Then, with $(\delta_1,\delta_2)$ defined in \eqref{eq.delta}, one finds for iteration $(k+1) \in \N{}$ that
  \begin{align*}
    \E[f(x_{k+1})] - \finf
      &\leq (1 - \delta_1\alpha_k) (\E[f(x_k)] - \finf) + \delta_2 \alpha_k^2 \\
      &=    \(1 - \frac{\delta_1a}{b+k}\)(\E[f(x_{k})] - \finf) + \frac{\delta_2a^2}{(b+k)^2}\\
      &\leq \(1 - \frac{\delta_1a}{b+k}\) \frac{\phi}{b+k} + \frac{\delta_2a^2}{(b+k)^2} \\
      &=    \frac{(b+k)\phi}{(b+k)^2} - \frac{\delta_1a \phi}{(b+k)^2} + \frac{\delta_2 a^2}{(b+k)^2}\\
      &=    \frac{(b+k-1)\phi}{(b+k)^2} - \frac{(\delta_1a-1) \phi}{(b+k)^2} + \frac{\delta_2a^2}{(b+k)^2}\\
      &\leq \frac{(b+k-1)\phi}{(b+k)^2} \leq \frac{\phi}{b+k+1},
  \end{align*}
  where the last equation follow from the definition of $\phi$ in \eqref{eq.phi} and the last inequality follows from the fact that $(z - 1)(z + 1) \leq z^2$ for any $z \in \R{}$. \ifthenelse{\technicalreport=0}{\Halmos}{}
\endproof

\ref{alg.trish} can also yield sublinear decrease of the expected optimality gap with fixed parameters.  However, this can only be guaranteed with the stronger assumption on the stochastic gradient estimates stipulated in Assumption~\ref{ass.nonconvex.diminishing} (specifically in \eqref{eq.nonconvex.diminishing.variance}).

\btheorem\label{th.PL.sublinear_2}
  Under Assumptions~\ref{ass.f}, \ref{ass.gH}, \ref{ass.nonconvex.diminishing}, and \ref{ass.PL}, if the pair $(a,b) \in \R{}_{>0} \times \R{}_{>0}$ is chosen such that $\alpha_k \leq \frac{4}{\gamma_2c}$ for all $k \in \N{}$, then for all $k \in \N{}$ the expected optimality gap satisfies
  \bequation\label{eq.sublinear_2}
    \E[f(x_k)] - \finf \leq \frac{\phi}{b+k},
  \eequation
  where
  \bequationNN
    \phi := \max\left\{(b+1) (f(x_1)-\finf), \frac{\delta_2 a^2}{\delta_1a-1}\right\} \in (0,\infty),
  \eequationNN
  with
  \bequationNN
    \delta_1 := \frac14 \gamma_2 c \in \(0,\frac1\alpha\right]\ \ \text{and}\ \ \delta_2 = \(\frac{\gamma_1^2}{\gamma_2^2} - \frac18\gamma_2\)M_g \in (0,\infty).
  \eequationNN
\etheorem
\proof{Proof.}
  As in the proof of Theorem~\ref{th.nonconvex.diminishing} (see \eqref{eq.nonconvex.diminishing.reduction}), it follows for all $k \in \N{}$ that
  \bequationNN
    \E_k[f(x_{k+1})] \leq f(x_k) - \frac18 \gamma_2 \alpha_k \|\nabla f(x_k)\|_2^2 + \(\frac{\gamma_1^2}{\gamma_2} - \frac18\gamma_2\) M_g \alpha_k^2,
  \eequationNN
  Noting that this inequality has the same form as that in \eqref{eq.same_as_next}, the remainder of the proof follows in the same manner as that for Theorem~\ref{th.PL.sublinear_1}. \ifthenelse{\technicalreport=0}{\Halmos}{}
\endproof

Finally in this section, let us consider the behavior of the algorithm under the following stronger assumption, which requires that the variance of the stochastic gradient estimates vanishes at a geometric rate.  Specifically, consider the following assumption.

\bassumption\label{ass.PL.geometric}
  The variances of the stochastic gradient estimates decreases at a geometric rate and the sequence of stochastic Hessian estimates is uniformly bounded in the sense that there exist constants $(M_g,M_H,\zeta) \in \R{}_{>0} \times \R{}_{>0} \times (0,1)$ such that, for all $k \in \N{}$,
  \bequationNN
    \E_k[\|\nabla f(x_k) - g_k\|_2^2] \leq M_g \zeta^{k-1}\ \ \text{and}\ \ \|H_k\|_2 \leq M_H.
  \eequationNN
  In addition, $(\alpha_k,\gamma_{1,k},\gamma_{2,k}) = (\alpha,\gamma_1,\gamma_2)$ for all $k \in \N{}$, where $\gamma_1 \geq \gamma_2 > 0$ and
  \bequationNN
    0 < \alpha \leq \frac{\gamma_2}{4\gamma_1^2(L_g + M_H)},
  \eequationNN
  which, in particular, implies that \eqref{eq.alpha} holds for all $k \in \N{}$.
\eassumption

This assumption leads to the following theorem.

\btheorem
  Under Assumptions \ref{ass.f}, \ref{ass.gH}, \ref{ass.PL}, and \ref{ass.PL.geometric}, \ref{alg.trish} yields
  \bequation\label{eq.linear_rate}
    \E[f(x_k)] - \finf \leq \omega \rho^{k-1},
  \eequation
  where
  \bequation\label{eq.kappa}
    \baligned
      \kappa_1 := \frac18\gamma_2,\ \ \kappa_2 := \(\frac{\gamma_1^2}{\gamma_2} - \frac18 \gamma_2\) M_g,\ \ \omega &:= \max\left\{f(x_1) - \finf,\frac{\kappa_2}{c\kappa_1}\right\}, \\
      \text{and}\ \ \rho &:= \max\{1 - c \kappa_1 \alpha, \zeta\} \in (0,1).
    \ealigned
  \eequation
\etheorem
\proof{Proof.}
  Using the same arguments as in the beginning of the proof of Theorem~\ref{th.nonconvex.fixed} (specifically leading to \eqref{eq.for_later}), one has for all $k \in \N{}$ that
  \bequationNN
    \E_k[f(x_{k+1})] \leq f(x_k) - \frac18 \gamma_2 \alpha \|\nabla f(x_k)\|_2^2 + \alpha \(\frac{\gamma_1^2}{\gamma_2} - \frac18 \gamma_2\) M_g \zeta^{k-1}.
  \eequationNN
  Applying the bound in Assumption~\ref{ass.PL}, subtracting $\finf$ from both sides, and taking total expectation, one finds with $(\kappa_1,\kappa_2)$ defined in \eqref{eq.kappa} that, for all $k \in \N{}$, one has
  \bequationNN
    \E[f(x_{k+1})] - \finf \leq (1 - 2c \kappa_1 \alpha) (\E[f(x_k)] - \finf) + \kappa_2 \alpha \zeta^{k-1}.
  \eequationNN
  Let us now prove \eqref{eq.linear_rate} by induction. First, for $k = 1$, the inequality follows by the definition of~$\omega$ in \eqref{eq.kappa}.  Then, assuming the inequality holds true for $k \in \N{}$, one finds from above that
  \bequationNN
    \baligned
      \E[f(x_{k+1})] - \finf
        &\leq (1 - 2 c \kappa_1 \alpha) \omega \rho^{k-1} + \kappa_2 \alpha \zeta^{k-1} \\
       &=    \omega \rho^{k-1} \(1 - 2 c \kappa_1 \alpha + \frac{\kappa_2 \alpha}{\omega} \(\frac{\zeta}{\rho}\)^{k-1}\) \\
       &\leq \omega \rho^{k-1} \(1 - 2 c \kappa_1 \alpha + \frac{\kappa_2 \alpha}{\omega}\) \leq \omega \rho^{k-1} (1 - c \kappa_1 \alpha) \leq \omega \rho^k,
    \ealigned
  \eequationNN
  which proves that the conclusion holds for $k + 1$, as desired. \ifthenelse{\technicalreport=0}{\Halmos}{}
\endproof

\section{Complexity Analysis}

In this section, we prove a complexity result for \ref{alg.trish}.  While not representing the behavior of the algorithm in the fully stochastic regime, the result does show that if one computes sufficiently accurate gradient and Hessian estimates, then one obtains---with the same algorithm---a worst-case performance that is reminiscent of results that can be proved for certain deterministic algorithms with optimal complexity properties.  To keep our result in the stochastic setting, we assume only that the stochastic gradients and Hessians are sufficiently accurate in expectation.  Consequently, our theorem is weaker than those that can be proved in the deterministic setting.  (If one were to replace the conditional expectations in \eqref{eq.complexity_derivatives} with computed values, then the same arguments would show that \ref{alg.trish} yields first-order $\epsilon$-stationarity in at most $\Ocal(\epsilon^{-3/2})$ iterations.)

\bassumption\label{ass.complexity}
  The Hessian function $\nabla^2 f : \R{n} \to \R{n \times n}$ is Lipschitz continuous with constant $L_H \in \R{}_{>0}$.  In addition, given $\epsilon \in \R{}_{>0}$, the expected distances of the stochastic gradient and stochastic Hessian estimates from the true gradients and Hessians, respectively, are uniformly bounded with respect to $(L_H,\epsilon)$ in the sense that there exist constants $\mu_1 \in (0,\tfrac{1}{12})$ and $\mu_2 \in (0,\tfrac{1}{12})$ such that, for all $k \in \N{}$,
  \bequation\label{eq.complexity_derivatives}
    \E_k[\|\nabla f(x_k) - g_k\|_2] \leq \frac{\mu_1}{L_H}\epsilon\ \ \text{and}\ \ \E_k[\|\nabla^2 f(x_k) - H_k\|_2] \leq \mu_2\sqrt{\epsilon}.
  \eequation
  Moreover, for all $k \in \N{}$, the subproblem~\eqref{prob.sub} is solved to global optimality.  Finally, the norms of the stochastic gradients are uniformly bounded above and below in that there exists $(G_{low},G_{high}) \in \R{}_{>0} \times \R{}_{>0}$ such that $G_{low} \leq \|g_k\|_2 \leq G_{high}$ for all $k \in \N{}$.
\eassumption 

Under Assumption~\ref{ass.complexity}, since the subproblem~\eqref{prob.sub} is solved to global optimality for all $k \in \N{}$, it follows for all $k \in \N{}$ that there exists a scalar $\upsilon_k$ such that
\bsubequations\label{prob.sub_optimality}
  \begin{align}
    g_k + (H_k + \upsilon_k I) s_k &= 0 \label{eq.complexity.dualopt} \\
    H_k + \upsilon_k I &\succeq 0 \label{eq.complexity.psd} \\
    \text{and}\ \ 0 \leq \upsilon_k \perp \Delta_k - \|s_k\|_2 &\geq 0. \label{eq.complexity.comp}
  \end{align}
\esubequations

\btheorem
  Suppose Assumptions~\ref{ass.f}, \ref{ass.gH}, and \ref{ass.complexity} hold.  In addition, suppose that $(\alpha_k,\gamma_{1,k},\gamma_{2,k}) = (\alpha,\gamma_1,\gamma_2)$ for all $k \in \N{}$, where $\gamma_1 \geq \gamma_2 > 0$ and for some constants $(\lambda_1,\lambda_2,\lambda_3) \in (0,1) \times (0,1) \times (0,1)$ satisfying
  \bequation\label{eq.complexity_parameter_choice}
    \lambda_1^2\lambda_2^2 - \frac{\mu_1}{\lambda_3} - \frac{\mu_2}{\lambda_3^2} - \frac{2}{3\lambda_3^3} \geq \frac16
  \eequation
  one has that the parameters employed by \ref{alg.trish} satisfy
  \bequation\label{eq.complexity.parameters}
    \alpha \in \left[\frac{2\lambda_1\sqrt{\epsilon}}{L_H},\frac{2\sqrt{\epsilon}}{L_H}\right],\ \ \gamma_1 \in \left[\frac{\lambda_2}{G_{low}},\infty\),\ \ \text{and}\ \ \gamma_2 \in \(0,\frac{1}{\lambda_3G_{high}}\right].
  \eequation
  Then, either $\upsilon_k > \sqrt{\epsilon}$ for all $k \in \{1,\dots,K\}$ where $K = \Ocal(\epsilon^{-3/2})$ and the (conditionally) expected total decrease in $f$ in these iterations is at least the initial optimality gap, i.e.,
  \bequationNN
    \sum_{k=1}^K \E_k[f(x_k) - f(x_{k+1}) | \upsilon_k > \sqrt{\epsilon}] \geq f_0 - \finf
  \eequationNN
  or for some $\overline{K} = \Ocal(\epsilon^{-3/2})$ one finds that
  \bequationNN
    \E_{\overline{K}}[\|\nabla f(x_{\overline{K}+1})\|_2 | \upsilon_{\overline{K}} \leq \sqrt{\epsilon}] \leq \Ocal(\epsilon).
  \eequationNN
\etheorem
\proof{Proof.}
  Under Assumption~\ref{ass.complexity}, it follows for all $k \in \N{}$ that
  \bsubequations
    \begin{align}
      f(x_k+ s_k) - f(x_k) - \nabla f(x_k)^T s_k - \frac{1}{2} s_k^T \nabla^2 f(x_k) s_k &\leq \frac{L_H}{6} \|s_k\|^3 \label{eq.complexity.func} \\
      \text{and}\ \ \|\nabla f(x_k + s_k) - \nabla f(x_k) - \nabla^2 f(x_k) s_k\| &\leq \frac{L_H}{2} \|s_k\|^2. \label{eq.complexity.grad}
    \end{align}
  \esubequations
  
  For the next parts of the proof, we consider two cases.  In the first case, we show that if the nonnegative scalar $\upsilon_k$ in the optimality conditions \eqref{prob.sub_optimality} is sufficiently small for some $k \in \N{}$, then the conditional expectation of the gradient of $f$ at $x_{k+1}$ is at most proportional to $\epsilon$.  In the second case, when $\upsilon_k$ is not sufficiently small, we show that the conditional expected decrease in the objective function value is at least proportional to $\epsilon^{3/2}$.
  
  First, suppose that $\upsilon_k \leq \sqrt{\epsilon}$.  For all $k \in \N{}$, it follows from the Cauchy-Schwarz inequality, \eqref{eq.complexity.grad}, \eqref{eq.complexity.dualopt}, and the trust region constraint in \eqref{prob.sub} that
  \begin{align}
    \|\nabla f(x_{k+1})\|_2
      &\leq \|\nabla f(x_{k+1}) - \nabla f(x_k) - \nabla^2 f(x_k) s_k\|_2 \nonumber \\
      &\qquad + \|\nabla f(x_k) - g_k\|_2 + \|(\nabla^2 f(x_k)-H_k) s_k\|_2 + \|g_k + H_k s_k\|_2 \nonumber \\
      &\leq \frac{L_H}{2} \|s_k\|_2^2 + \|\nabla f(x_k) - g_k\|_2 + \|\nabla^2 f(x_k)-H_k\|_2 \|s_k\|_2 + \upsilon_k \|s_k\|_2 \nonumber \\
      &\leq \frac{L_H}{2} \Delta_k^2 + \|\nabla f(x_k) - g_k\|_2 + \Delta_k \|\nabla^2 f(x_k)-H_k\|_2 + \Delta_k\sqrt{\epsilon}.\label{eq.complexity.newgrad}
  \end{align}
  Let us now consider the three cases defined on page~\pageref{case.1}.  In \ref{case.1}, one has that $\|g_k\|_2 \leq 1/\gamma_1$, meaning $\Delta_k = \gamma_1 \alpha \|g_k\|_2 \leq \alpha$.  In \ref{case.2}, one has that $\Delta_k = \alpha$.  Finally, in \ref{case.3}, one has by \eqref{eq.complexity.parameters} that $\Delta_k = \gamma_2 \alpha \|g_k\|_2 \leq \gamma_2 \alpha G_{high} \leq \alpha/\lambda_3$.  Thus, it follows by \eqref{eq.complexity.newgrad} that, for all $k \in \N{}$,
  \bequationNN
    \baligned
      \|\nabla f(x_{k+1})\|_2 &\leq \frac{L_H}{2} \(\frac{\alpha}{\lambda_3}\)^2 + \|\nabla f(x_k) - g_k\|_2 + \frac{\alpha}{\lambda_3} \|\nabla^2 f(x_k) - H_k\|_2 + \frac{\alpha}{\lambda_3} \sqrt{\epsilon} \nonumber \\
      &\leq \(\frac{1}{\lambda_3^2} + \frac{1}{\lambda_3}\) \frac{2}{L_H} \epsilon + \|\nabla f(x_k) - g_k\|_2 + \frac{2}{L_H\lambda_3} \|\nabla^2 f(x_k)-H_k\|_2 \sqrt{\epsilon}.
    \ealigned
  \eequationNN
  Taking conditional expectation, it follows for all $k \in \N{}$ that
  \bequationNN
    \baligned
      &\ \E_k[\|\nabla f(x_{k+1})\|_2 | \upsilon_k \leq \sqrt{\epsilon}] \\
      \leq&\ \(\frac{1}{\lambda_3^2} + \frac{1}{\lambda_3}\) \frac{2}{L_H} \epsilon + \E_k[\|\nabla f(x_k) - g_k\|_2] + \frac{2}{L_H\lambda_3} \E_k[\|\nabla^2 f(x_k)-H_k\|_2] \sqrt{\epsilon} \\
      \leq&\ \(\(\frac{1}{\lambda_3^2} + \frac{1}{\lambda_3}\) \frac{2}{L_H} + \frac{\mu_1}{L_H} + \frac{2\mu_2}{L_H\lambda_3}\) \epsilon.
    \ealigned
  \eequationNN
  
  Second, suppose that $\upsilon_k > \sqrt{\epsilon}$.  For such $k \in \N{}$, it follows by \eqref{eq.complexity.comp} that $\|s_k\|_2 = \Delta$.  Therefore, by \eqref{eq.complexity.func}, \eqref{prob.sub_optimality}, and the Cauchy-Schwarz inequality, it follows for all $k \in \N{}$ that
  \begin{align}
        &\ f(x_{k+1}) - f(x_k) \nonumber \\
    \leq&\ \nabla f(x_k)^T s_k + \frac{1}{2} s_k^T \nabla^2 f(x_k) s_k + \frac{L_H}{6} \|s_k\|_2^3 \nonumber \\
    \leq&\ g_k^T s_k + \half s_k^T H_k s_k + (\nabla f(x_k) - g_k)^T s_k + \half s_k^T(\nabla^2 f(x_k) - H_k)s_k + \frac{L_H}{6} \|s_k\|_2^3 \nonumber \\
    \leq&\ g_k^T s_k + \half s_k^T H_k s_k + \|\nabla f(x_k) - g_k\|_2 \|s_k\|_2 + \half \|\nabla^2 f(x_k) - H_k\|_2 \|s_k\|_2^2 + \frac{L_H}{6} \|s_k\|_2^3 \nonumber \\
       =&\ g_k^T s_k + \half s_k^T H_k s_k +  \Delta_k \|\nabla f(x_k) - g_k\|_2 + \frac{1}{2} \Delta_k^2 \|\nabla^2 f(x_k) - H_k\|_2 + \frac{L_H}{6} \Delta_k^3 \nonumber \\
    \leq&\ -\half \upsilon_k \Delta_k^2 + \Delta_k \|\nabla f(x_k) - g_k\|_2 + \half \Delta_k^2 \|\nabla^2 f(x_k) - H_k\|_2 + \frac{L_H}{6} \Delta_k^3 \nonumber \\
    \leq&\ -\half \sqrt{\epsilon} \Delta_k^2 + \Delta_k \|\nabla f(x_k) - g_k\|_2 + \half \Delta_k^2 \|\nabla^2 f(x_k) - H_k\|_2 + \frac{L_H}{6} \Delta_k^3. \label{eq.complexity.decrease}
  \end{align}
  Let us consider the three cases defined on page~\pageref{case.1}.  In \ref{case.1}, it follows under Assumption~\ref{ass.complexity} and by \eqref{eq.complexity.parameters} and the fact that $\|g_k\|_2 \leq 1/\gamma$ that $\Delta_k = \gamma_1 \alpha \|g_k\|_2 \geq \gamma_1 \alpha G_{low} \geq \frac{2\lambda_1\lambda_2}{L_H} \sqrt{\epsilon}$ and $\Delta_k \leq \alpha \leq \frac{2}{L_H} \sqrt{\epsilon}$.  In \ref{case.2}, one finds that $\Delta_k = \alpha \in \left[\frac{2\lambda_1}{L_H}\sqrt{\epsilon},\frac{2}{L_H}\sqrt{\epsilon}\right]$.  Finally, in \ref{case.3}, one finds as before that $\Delta_k \leq \alpha/\lambda_3$, meaning that $\Delta_k \leq \frac{2}{L_H\lambda_3} \sqrt{\epsilon}$.  Moreover, one finds by the fact that $\|g_k\|_2 \geq 1/\gamma_2$ in this case that $\Delta_k \geq \alpha \geq \frac{2\lambda_1}{L_H}\sqrt{\epsilon}$.  Hence, for all $k \in \N{}$, one finds
  \bequationNN
    \Delta_k \in \left[\frac{2\lambda_1\lambda_2}{L_H} \sqrt{\epsilon}, \frac{2}{L_H\lambda_3} \sqrt{\epsilon} \right].
  \eequationNN
  Combining this inclusion with \eqref{eq.complexity.decrease} and \eqref{eq.complexity_parameter_choice}, one finds for all $k \in \N{}$ that
  \bequationNN
    \baligned
      \E_k[f(x_{k+1}) | \upsilon_k > \sqrt{\epsilon}] - f(x_k)
      \leq&\ - \frac{2\lambda_1^2\lambda_2^2}{L_H^2} \epsilon^{3/2} + \frac{2\mu_1}{L_H^2\lambda_3} \epsilon^{3/2} + \frac{2\mu_2}{L_H^2\lambda_3^2} \epsilon^{3/2} + \frac{4}{3L_H^2\lambda_3^3} \epsilon^{3/2} \\
      \leq&\ -\frac{2}{L_H^2}\(\lambda_1^2\lambda_2^2 - \frac{\mu_1}{\lambda_3} - \frac{\mu_2}{\lambda_3^2} - \frac{2}{3\lambda_3^3}\) \epsilon^{3/2} \\
      \leq&\ -\frac{1}{3L_H^2} \epsilon^{3/2}.
    \ealigned
  \eequationNN
  
  Combining the results of these two cases leads to our desired conclusion.  In particular, suppose that $\upsilon_k > \sqrt{\epsilon}$ for all $k \in \{1,\dots,\khat\}$ for some $\khat \in \N{}$.  It follows from above that
  \bequationNN
    \sum_{k=1}^{\khat} \E_k[f(x_k) - f(x_{k+1}) | \upsilon_k > \sqrt{\epsilon}] \geq \khat \(\frac{1}{3L_H^2}\) \epsilon^{3/2}.
  \eequationNN
  The left-hand side of this inequality is greater than $f_1 - \finf$ as long as
  \bequationNN
    \khat \(\frac{1}{3L_H^2}\) \epsilon^{3/2} \geq f_1 - \finf \iff \khat \geq 3L_H^2(f_1 - \finf) \epsilon^{-3/2},
  \eequationNN
  which shows that the desired inequality is true for some $K = \Ocal(\epsilon^{-3/2})$. \ifthenelse{\technicalreport=0}{\Halmos}{}
\endproof

\section{Numerical Experiments}

The goal of our numerical experiments is to show that TRish, with stochastic second-order derivative information incorporated, can outperform SG and first-order TRish (i.e., TRish with $H_k=0$ for all $k \in \N{}$).  In particular, our goal is to show with a few interesting test problems that TRish can offer a better final solution, and offer better stability throughout the optimization process, in the sense that the quality of the solution estimates does not vary as wildly from one iteration to the next as it might for SG.

\subsection{Implementation Details}

We implemented TRish and SG in Python.  All of our test problems involve training neural networks.  The problems were implemented using PyTorch, which allows one to use back propagation to compute stochastic gradient estimates and perform matrix-vector products with stochastic Hessian estimates.  For TRish, we implemented a Steihaug-CG routine (see \cite{Stei83}) for approximately solving the trust region subproblems, where for each subproblem the same batch of data samples used to define the stochastic gradient estimate is used to define the stochastic Hessian estimate.  To ensure that TRish did not expend too much effort solving any single subproblem, we imposed a limit of 3 on the number of CG iterations performed when solving each subproblem.  In our comparisons, we equate the cost of one stochastic gradient estimate with the cost of computing one stochastic-Hessian-vector product.  This allows SG and first-order TRish to perform more optimization iterations per epoch than TRish is able to perform.

\subsection{Hyperparameters Tuning}

The hyperparameters for all algorithms were tuned using a similar approach to that used in~\cite{CurtScheShi19}.  In particular, for each test problem, we proceeded as follows.  First, to establish a baseline for the hyperparameter values, we ran SG with a fixed stepsize of $\alpha = 0.1$ and computed $G$ as the average norm of the stochastic gradient estimates computed throughout the run.  We then established sets of possible hyperparameter values with the formulas $\alpha = 10^{\lambda}$, $\gamma_1 = \frac{2^a}{G}$, and $\gamma_2=\frac{1}{2^bG}$, where $\lambda$, $a$, and $b$ were evenly distributed in some interval.  (Different intervals were used for each test problem so that, e.g., the best stepsize for SG was never at the extreme of the allowed range.  Details are given in the following subsections for each test problem.)  For simplicity, we only consider the behavior of the algorithms with fixed hyperparameter values.  To ensure that all algorithms were tuned with the same amount of effort, we fixed the total number of hyperparameter settings to be the same for all algorithms.  For example, if (first-order) TRish considers 4 values of $\alpha$, 3 values of $\gamma_1$, and 3 values of $\gamma_2$, then we allowed SG to consider $4 \times 3 \times 3 = 36$ stepsizes.

To choose the best hyperparameter values for each algorithm for each test problem, we used a standard type of cross validation procedure.  Each dataset came equipped with a training set and a testing set of data.  We began by randomly selecting points from the training set to form a validation set.  For each hyperparameter setting, we ran each algorithm and observed its performance in terms of final validation accuracy (in the case of image classification) or final validation loss (in terms of time series forecasting).  Once the best hyperparameter setting was found in this manner, we ran the algorithm using this setting on \emph{all} of the original training data.  In the subsections below, we provide plots of the accuracy and/or loss during this final run for the training and testing data.

\subsection{\texttt{FashionMNIST}}

The first dataset that we considered was \texttt{FashionMNIST} (\cite{xiao2017fashion}).  This consists of images of 10 different types of clothing.  Each image is a color image of size $28 \times 28$.  There are 60000 training images and 10000 testing images.  We randomly chose 10000 images out of the training set as our validation set, and chose the best set of hyperparameters for each algorithm as the one yielding highest classification accuracy on the validation set.

The neural network that we considered for performing classification for this dataset was composed of two convolutional layers (involving 10 and 20 output channels, respectively, with kernal size 5) followed by a dropout layer and three fully connected layers.  ReLU activation was used at each hidden layer and the objective is defined using the logistic loss (cross entropy) function.  It is known that one can achieve better classification accuracy on \texttt{FashionMNIST} using a more sophisticated neural network, but this network offers sufficiently good results in order for us to demonstrate the behavior of TRish.

We ran each algorithm for 5 epochs with a mini-batch size of 128.  During tuning, we obtained $G=1.5644$. For TRish and first-order TRish, we considered 8 stepsize values over $[0.1,1]$, namely, $\alpha = 10^{-1+i/7}$ for $i \in \{0,1,\dots,7\}$, along with $\gamma_1 \in \{\tfrac4G,\tfrac{16}G\}=\{2.5568,10.2274\}$ and $\gamma_2 \in \{\tfrac1{2G},\tfrac1{8G}\} = \{0.3196,0.07990\}$.  For a fair comparison (see \cite{CurtScheShi19}), this means that it was appropriate to allow SG to consider 32 stepsize choices in the range $[\tfrac1{8G}\times10^{-1},\tfrac{16}G\times10^0] = [10^{-2.0974},10^{1.0097}] = [0.00799, 10.2275]$.  TRish ended up with the values $(\alpha,\gamma_1,\gamma_2) = (0.1930,10.2274,0.07990)=(10^{-5/7},\frac{16}{1.5644},\frac{1}{8(1.5644)})$, first-order TRish ended up with the values $(\alpha,\gamma_1,\gamma_2) = (0.3727,2.5568,0.3196)=(10^{-3/7},\frac{4}{1.5644},\frac{1}{2(1.5644)})$, and SG ended up with the value $\alpha = 0.4192=10^{-0.3775}$.

Once the hyperparameter values were determined, we ran the algorithms on the training data 5 times each.  In Figure~\ref{fig.FashionMNIST}, we plot the training loss and testing accuracy over the 5 epochs.  The line for each algorithm for each plot shows the mean values over the 5 runs with the shaded region showing one standard deviation above and below the mean.  One finds that while the first-order algorithms have an edge in the early parts of the runs, eventually TRish overtakes both of the other algorithms in terms of final training loss (for which lower is better) and final testing accuracy (for which higher is better).

\bfigure[ht]
  \centering
  \begin{subfigure}{.49\textwidth}
    \centering
    \includegraphics[width=1\linewidth]{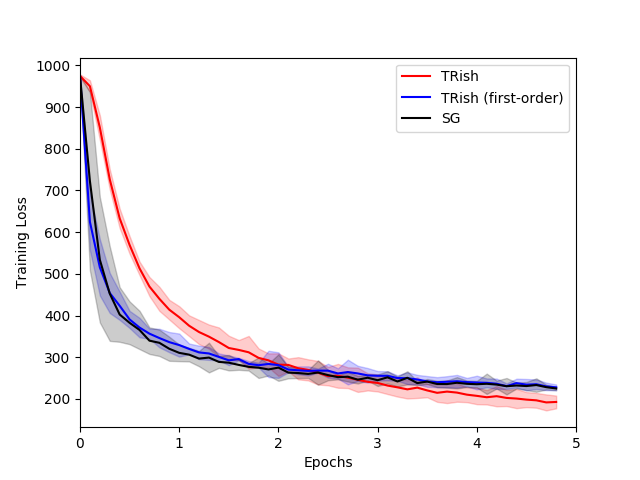}
  \end{subfigure}
  \begin{subfigure}{.49\textwidth}
    \centering
    \includegraphics[width=1\linewidth]{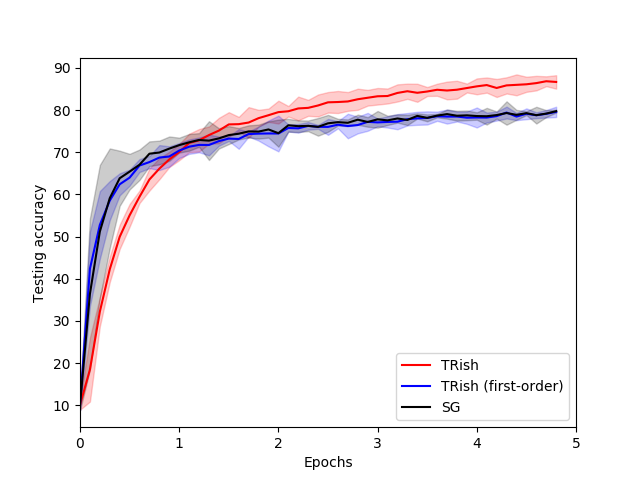}
  \end{subfigure}
  \caption{Training loss and testing accuracy during the first five epochs when \ref{alg.trish}, first-order TRish, and SG are employed to train a convolutional neural network over the \texttt{FashionMNIST} dataset.}
  \label{fig.FashionMNIST}
\efigure

\subsection{\texttt{CIFAR-10}.}

The second dataset that we considered was \texttt{CIFAR-10} (\cite{CIFAR10}).  This dataset consists of 10 classes of color images of different objects.  Each image has size $32 \times 32$.  There are 50000 training images and 10000 testing images.  We randomly chose 5000 of the training images to compose the validation set.  As in the previously subsection, the best set of hyperparameters for each algorithm was chosen as the one yielding highest classification accuracy on the validation set.

The neural network that we considered for this dataset was composed of two convolutional layers (involving 6 and 16 output channels, respectively, with kernal size 5) followed by a max pooling layer, a dropout layer, and three fully connected layers.  ReLU activation was used at each hidden layer and the objective was again the logistic loss function.  Again, one can achieve better testing accuracy using a more sophisticated neural network, but this network gave sufficiently good results to demonstrate the behavior of our algorithm.

We ran 5 epochs with a mini-batch size of 128.  We obtained $G=2.7819$ and considered $\alpha = 10^{-1+i/7}$ for $i \in \{0,1,\dots,7\}$, $\gamma_1 \in \{\tfrac4G,\tfrac{16}G\}=\{1.4378,5.7515\}$, and $\gamma_2 \in \{\tfrac1{4G},\tfrac1{80G}\}=\{0.08986,0.004493\}$.  This means that SG was tuned with 32 choices of $\alpha$ in the range $[\tfrac1{80G}\times10^{-1},\tfrac{16}G\times10^0]=[10^{-3.3474},10^{0.7598}]=[0.0004493, 5.7515]$.  TRish chose $(\alpha,\gamma_1,\gamma_2) = (0.1389,5.7515,0.004493)=(10^{-6/7},\frac{16}{2.7819},\frac{1}{80(2.7819)})$, first-order TRish chose $(\alpha,\gamma_1,\gamma_2) = (0.3727,5.7515,0.08986)=(10^{-3/7},\frac{16}{2.7819},\frac{1}{4(2.7819)})$, and SG chose $\alpha = 0.2316=10^{-0.6352}$.

Figure~\ref{fig.CIFAR10} shows the result of this experiment over 5 runs.  Interestingly, for this problem, TRish does not outperform the others in terms of training loss; indeed, first-order TRish appears to give the best results in terms of training loss.  However, TRish eventually offers better testing accuracy.  While one cannot guarantee that such would be the behavior in general, one does see benefits of TRish-based methods compared to SG.

\bfigure[ht]
  \centering
  \begin{subfigure}{.49\textwidth}
    \centering
    \includegraphics[width=1\linewidth]{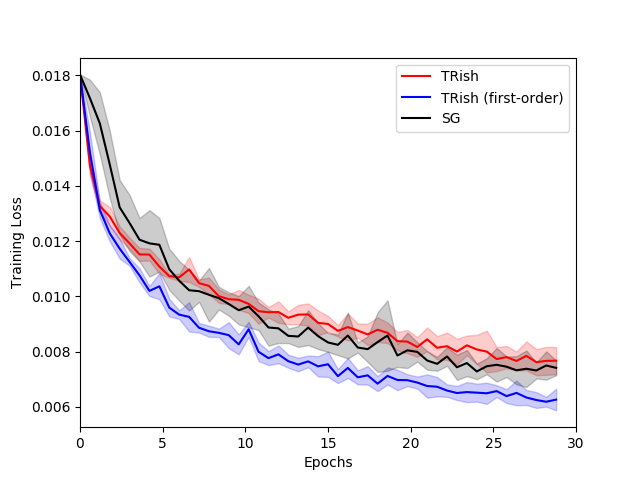}
  \end{subfigure}
  \begin{subfigure}{.49\textwidth}
    \centering
    \includegraphics[width=1\linewidth]{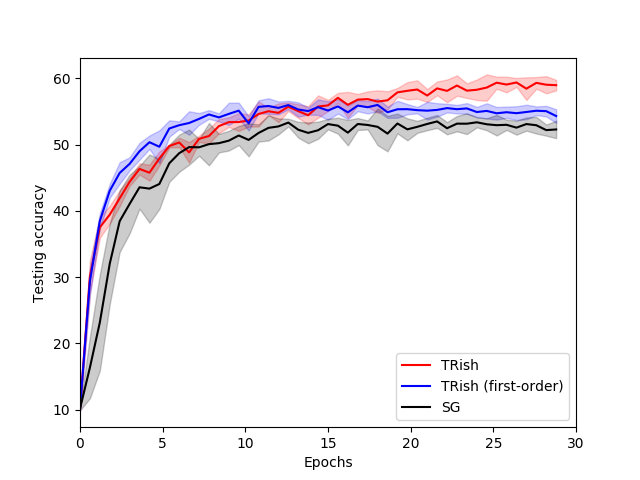}
  \end{subfigure}
  \caption{Training loss and testing accuracy during the first five epochs when \ref{alg.trish}, first-order TRish and SG are employed to train a convolutional neural network over the \texttt{CIFAR10} dataset.}
  \label{fig.CIFAR10}
\efigure

\subsection{NSW2016}

As a final test problem, we considered one of time series forecasting.  For this, we used historical data posted online by the Australian Energy Market Operator (AEMO) on demand for electricity in New South Wales in 2016.\footnote{\href{https://www.aemo.com.au/Electricity/National-Electricity-Market-NEM/Data-dashboard}{https://www.aemo.com.au/Electricity/National-Electricity-Market-NEM/Data-dashboard}}  This gives a univariate time series of length 17423.  We used the first 17000 values for our experiments.  We used the first 12000 as the training set, the following 2000 as the validation set, and the remaining 3000 as the testing set.  We chose the set of hyperparameters that yielded the lowest validation loss.

The recurrent neural network that we considered for this dataset was composed of a single long short-term memory (LSTM) layer with hidden size 32 followed by a fully connected layer.  A time step of 10 was used with ReLU activation after the LSTM layer.  The objective function used was the mean squared error.

We ran the experiment for 20 epochs using a mini-batch size of 100.  We obtained $G=720.1389$ and considered $\alpha = 10^{-1+i/3}$ for $i = \{0,1,\dots,6\}$ along with $\gamma_1 \in \{\tfrac4G,\tfrac{16}G\}=\{0.005555,0.02222\}$ and $\gamma_2 \in \{\tfrac1{2G},\tfrac1{20G}\}=\{0.0006944,0.00006944\}$.  SG was tuned with 16 choices of $\alpha$ in the range $[\tfrac1{20G}\times10^{-1},\tfrac{16}G\times10^1]=[0.000006944, 0.2222]=[10^{-5.1586},10^{-0.6532}]$.  As a result of hyperparameter tuning, TRish chose $(\alpha,\gamma_1,\gamma_2) = (2.1544,0.2222,0.00006944)=(10^{1/3},\frac{16}{720.1389},\frac{1}{20(720.1389)})$, first-order TRish chose $(\alpha,\gamma_1,\gamma_2) = (0.4641,0.005555,0.0006944)=(10^{-1/3},\frac{4}{720.1389},\frac{1}{2(720.1389)})$, and SG chose $\alpha = 0.0002204=10^{-3.6567}$.

Figure~\ref{fig.NSW2016} shows the result of this experiment over 50 runs.  The losses are plotted on a logarithmic scale for better viewing of the differences.  From the plots, it is clear that while all algorithms reach solutions of comparable quality, TRish is able to achieve low losses earlier than the other two methods, and first-order TRish similarly outperforms SG.

\bfigure[ht]
  \centering
  \begin{subfigure}{.49\textwidth}
    \centering
    \includegraphics[width=1\linewidth]{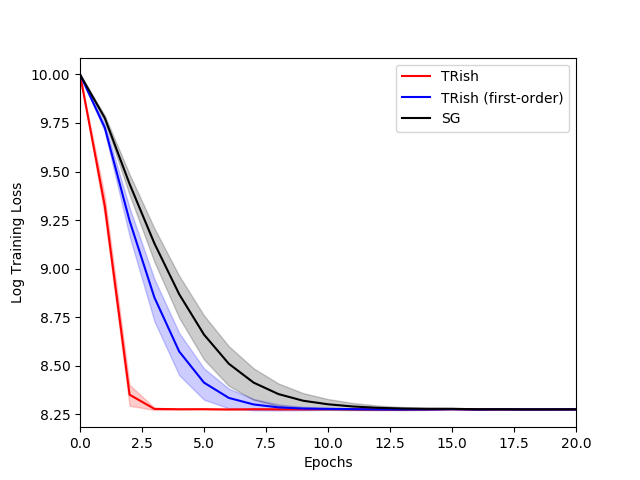}
  \end{subfigure}
  \begin{subfigure}{.49\textwidth}
    \centering
    \includegraphics[width=1\linewidth]{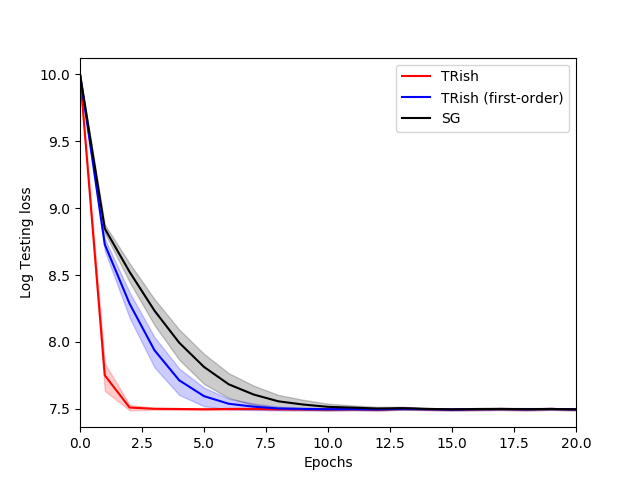}
  \end{subfigure}
  \caption{Training loss and testing loss during the first twenty epochs when \ref{alg.trish}, first-order TRish, and SG are employed to train a recurrent neural network over the \texttt{NSW2016} dataset.}
  \label{fig.NSW2016}
\efigure

\section{Conclusion}

A stochastic second-order trust region algorithm has been proposed, analyzed, and tested.  It can be viewed as a second-order extension of the algorithm proposed in~\cite{CurtScheShi19}.  We proved theoretical guarantees for the method that are on par with those proved for the first-order algorithm in~\cite{CurtScheShi19}, and in turn comparable to those possessed by SG and many of its variants.  That said, our numerical experiments demonstrate that the algorithm can perform better in practice, in terms of reaching better solutions and with more stable behavior.  We attribute this better behavior to the algorithm's use of carefully chosen trust region radii and stochastic second-order information.

\bibliographystyle{plain}
\bibliography{references}

\end{document}